\documentclass[12pt, a4paper]{article}
\title{On tensor products of matrix factorizations}
\author{Yves Baudelaire Fomatati\\
\small Department of Mathematics and Statistics, University of Ottawa,\\ \small Ottawa, Ontario, Canada K1N 6N5.\\ \small yfomatat@uottawa.ca.}

\date{}
\usepackage{latexsym}
\usepackage{amsmath,amssymb}
\usepackage[square,comma,numbers,sort&compress]{natbib}
\usepackage{geometry,mathtools}
\usepackage{resizegather}

\usepackage{leqno}
\usepackage{amsfonts}

\usepackage{enumerate}

\usepackage{amsthm}
\usepackage[utf8,latin9]{inputenc}

\newcommand{\rvline}{\hspace*{-\arraycolsep}\vline\hspace*{-\arraycolsep}}

\theoremstyle{plain}
\newtheorem{remark}{Remark}[section]
\theoremstyle{plain}
\newtheorem{lemma}{Lemma}[section]
\theoremstyle{plain}
\newtheorem{proposition}{Proposition}[section]
\theoremstyle{plain}
\newtheorem{theorem}{Theorem}[section]
\theoremstyle{plain}
\newtheorem{definition}{Definition}[section]
\theoremstyle{plain}
\newtheorem{corollary}{Corollary}[section]
\theoremstyle{plain}
\newtheorem{fact}{Fact}[section]
\theoremstyle{plain}

\newtheorem{example}{Example}[section]

\theoremstyle{plain}

\usepackage[all, 2cell]{xy}
\usepackage{txfonts}

\begin{document}
\maketitle
\begin{quote}
  \textbf{Abstract}
\end{quote}
Let $K$ be a field. Let $f\in K[[x_{1},...,x_{r}]]$ and $g\in K[[y_{1},...,y_{s}]]$ be nonzero elements. If $X$ (resp. $Y$) is a matrix factorization of $f$ (resp. $g$), Yoshino had constructed a tensor product (of matrix factorizations) $\widehat{\otimes}$ such that $X\widehat{\otimes}Y$ is a matrix factorization of $f+g\in K[[x_{1},...,x_{r},y_{1},...,y_{s}]]$. In this paper, we propose a bifunctorial operation $\widetilde{\otimes}$ and its variant $\widetilde{\otimes}'$ such that $X\widetilde{\otimes}Y$ and $X\widetilde{\otimes}' Y$ are two different matrix factorizations of $fg\in K[[x_{1},...,x_{r},y_{1},...,y_{s}]]$. We call $\widetilde{\otimes}$ \textit{the multiplicative tensor product} of $X$ and $Y$. Several properties of $\widetilde{\otimes}$ are proved. Moreover, we find three functorial variants of Yoshino's tensor product $\widehat{\otimes}$. Then, $\widetilde{\otimes}$ (or its variant) is used in conjunction with $\widehat{\otimes}$ (or any of its variants) to give an improved version of the standard algorithm for factoring polynomials using matrices on the class of \textit{summand-reducible polynomials} defined in this paper. Our algorithm produces matrix factors whose size is at most one half the size one obtains using the standard method.
\section{Introduction}
In 1980, Eisenbud introduced the concept of matrix factorization in his seminal paper \cite{eisenbud1980homological}. He came up with an approach to factorize all polynomials including irreducible ones using matrices. For example the irreducible polynomial $f(x)=x^{2}+1$ over $\mathbb{R}[x]$ can be factored as follows: $$\begin{bmatrix}
    x  &  -1      \\
    1  &  x
\end{bmatrix}
\begin{bmatrix}
    x  &  1      \\
    -1  & x
\end{bmatrix}
= (x^{2} + 1)\begin{bmatrix}
    1  &  0      \\
    0  &  1
\end{bmatrix}
=fI_{2} $$
Thus; we say that
$
(\begin{bmatrix}
    x  &  -1      \\
    1  &  x
\end{bmatrix},
\begin{bmatrix}
    x  &  1      \\
    -1  & x
\end{bmatrix})
$
 is a $2\times 2$ matrix factorization of $f$.\\
 In a sense, this notion of factorizing polynomials using matrices can be seen as a generalization of the classical notion of polynomial factorization where a polynomial $h(x)=s(x)r(x)$ can be considered as the product of two $1\times 1$ matrices.\\
 There is a standard method for factoring (both reducible and irreducible) polynomials using matrices (cf. subsection \ref{sec: the std method}). One conspicuous downside of this algorithm is that for each monomial that is added to the polynomial, the size of the matrix factors doubles. As will be seen below (subsection \ref{sec: the std method}), polynomials with $n$ monomials have matrix factors of size $2^{n-1}$ (when $n=10$, we have matrix factors of size $2^{10-1}= 512$). Our main objective in this paper is to improve the standard method by creating a procedure that will yield matrix factors of smaller sizes than the ones one would normally obtain with the standard method.
 In \cite{crisler2016matrix}, the standard method was improved on a subclass of polynomials that are sums of squares. In this paper, we will improve the standard method on the class of \textit{summand-reducible polynomials} (cf. definition \ref{defn summand reducible polynomials}). Since a polynomial is made up of sums and products of monomials (or even other polynomials), one reasonable approach is to create two tensor products: The first with the ability to produce a matrix factorization of the sum of two polynomials from their respective matrix factorizations, and the second with the ability to produce a matrix factorization of the product of two polynomials from their respective matrix factorizations. The former was already constructed by Yoshino \cite{yoshino1998tensor} in 1998 as will be discussed below (cf. subsection \ref{subsec: Yoshino tens prod and variant}). So, one of our goals is to construct the latter.\\
 In the sequel, except otherwise stated, $K$ is a field and $K[[x]]$ (resp. $K[[y]]$) is the formal power series ring in the variables $x=x_{1},\cdots,x_{r}$ (resp. $y=y_{1},\cdots,y_{s}$). Let $f\in K[[x]]$ and $g\in K[[y]]$ be nonzero noninvertible\footnote{Yoshino \cite{yoshino1998tensor} requires an element $f\in K[[x]]$ to be noninvertible because if
$f$ is a unit, then $K[[x]]/(f)=$K[[x]]/K[[x]]$=\{1\}$. But in this paper we will not bother about such a restriction because we will not deal with the homological methods used in \cite{yoshino1998tensor}.} elements.\\
 Eisenbud also found out that matrix factorizations of a power series $f\in K[[x]]$ are closely related to the homological properties of modules over quotient rings $K[[x]]/(f)$. He proved that all maximal Cohen-Macaulay modules (MCM modules) without free summands are described by matrix factorizations.
 See \cite{leuschke2012cohen} and \cite{huneke2002two} for more background on MCM modules.
\\  Yoshino \cite{yoshino1998tensor} found a way to relate MCM modules over $K[[x]]/(f)$ and over $K[[y]]/(g)$ with MCM modules over $K[[x,y]]/(f+g)$. In fact, he constructed a tensor product denoted $\widehat{\otimes}$ which is such that if $X$ is a matrix factorization of $f\in K[[x]]$ and $Y$ is a matrix factorization of $g\in K[[y]]$, then $X\widehat{\otimes}Y$ is a matrix factorization of $f+g\in K[[x,y]]$. \\
In this paper, we are not interested in studying MCM modules and their relationship with some given quotient rings. Instead, without resorting to homological methods we construct a bifunctorial operation and its variant respectively denoted by $\widetilde{\otimes}$ and $\widetilde{\otimes}'$ which are such that
$X \widetilde{\otimes} Y$ and $X\widetilde{\otimes}'Y$ are two different (but of the same size) matrix factorizations of the product $fg\in K[[x,y]]$. \\

 Thus, our first main result is the following:
 \begin{quote}
   \textbf{Theorem A.}
1.
If $X$ (resp. $Y$) is a matrix factorization of $f$ (resp. $g$). Then, there is a tensor product $\widetilde{\otimes}$ of matrix factorizations which produces a matrix factorization $X \widetilde{\otimes} Y$ of the product $fg\in K[[x_{1},...,x_{r},y_{1},...,y_{s}]]$. $\widetilde{\otimes}$ is called the multiplicative tensor product of matrix factorizations. Moreover, $\widetilde{\otimes}$ has a functorial variant having the same effect on $X$ and $Y$.\\
2. The multiplicative tensor product $(-) \widetilde{\otimes} (-):MF(K[[x]],f)\times MF(K[[y]],g)\rightarrow MF(K[[x,y]],fg)$ is a bifunctor. Moreover, its variant $\widetilde{\otimes}'$ is also a bifunctor.
 \end{quote}

We also find three bifunctorial variants of Yoshino's tensor product $\widehat{\otimes}$. We then use the newly defined operation $\widetilde{\otimes}$ (or its variant) together with the existing Yoshino tensor product $\widehat{\otimes}$ (or any of its variants), to improve the standard method of matrix factorization of polynomials on the class of \textit{summand-reducible polynomials}, which is defined in this paper (cf. definition \ref{defn summand reducible polynomials}).\\

Our second main result is stated as follows:
 \begin{quote}
   \textbf{Theorem B.}
Let $f=t_{1}+\cdots + t_{s}+ g_{11}\cdots g_{1m_{1}} + \cdots + g_{l1}\cdots g_{lm_{l}}$ be a \textit{summand-reducible polynomial}. Let $p_{ji}$ be the number of monomials in $g_{ji}$. Then
there is an improved version of the standard method for factoring $f$ which produces factorizations of size $$2^{\prod_{i=1}^{m_{1}}p_{1i} + \cdots + \prod_{i=1}^{m_{l}}p_{li} - (\sum_{i=1}^{m_{1}}p_{1i} + \cdots + \sum_{i=1}^{m_{l}}p_{li})}$$ times smaller than the size one would normally obtain with the standard method.\\
 \end{quote}
Matrix factorizations is a rapidly growing field of research because it plays an important role in many areas of pure mathematics and physics.
In 1987, Buchweitz et al. \cite{buchweitz1987cohen} found that matrix factorizations of polynomials (over the reals) of the form $f_{n}=x_{1}^{2}+\cdots + x_{n}^{2}$, for $n=1,2,4$ and $8$ are related to the existence of composition algebras over $\mathbb{R}$ of dimension $1,2,4$ and $8$ (namely the complex numbers, the quaternions and the octonians).
It is a classical tool in the study of hypersurface singularity algebras (cf. \cite{eisenbud1980homological}). In 2002 and 2003,  Kapustin and Li in their papers \cite{kapustin2003topological} and \cite{kapustin2004d}, used matrix factorizations in string theory to study boundary conditions for strings in Landau-Ginzburg models. These models are very important in physics. The initial model was used to describe superconductivity \cite{patashinskii1973longitudinal}. Some of these models are also used in the field of research of mirror theory \cite{kapustin2004d}. A major advance was made by Orlov (\cite{orlov2003triangulated}, \cite{orlov2009derived}, \cite{orlov2012matrix}, \cite{orlov2006triangulated}), who showed that matrix factorizations could be used to study Kontsevich's homological mirror symmetry by giving a new description of singularity categories. Matrix factorizations have also proven useful for the study of cluster tilting \cite{dao2013vanishing}, Cohen-Macaulay modules and singularity theory (\cite{herzog1991linear}, \cite{buchweitz1987cohen}), Khovanov-Rozansky homology (\cite{khovanov2008matrix}), moduli of curves \cite{polishchuk2011matrix}, quiver and group representations (\cite{aspinwall2012quivers}, \cite{avramov1989modules}). In 2013, Yu \cite{yu2013geometric} in his PhD dissertation studied the geometry of the category of matrix factorizations. In 2014, Camacho \cite{camacho2015matrix} in chapter 4 of her PhD dissertation recalled the notion of graded matrix factorizations with special emphasis on $\mathbb{C}-$graded matrix factorizations.
Carqueville and Murfet in their paper \cite{carqueville2016adjunctions} published in 2016 briefly presented the construction of the bicategory $\mathcal{LG}_{K}$ of Landau-Ginzburg models whose $1$-cells are matrix factorizations. In the same year, another publication \cite{crisler2016matrix} from Crisler and Diveris examined matrix factorizations of polynomials in the ring $\mathbb{R}[x_{1},\cdots,x_{n}]$, using only techniques from elementary linear algebra. They focused mostly on factorizations of sums of squares of polynomials. They improved the standard method for factoring polynomials for a subclass of this class of polynomials.\\
The rest of this paper is organized as follows: In section $2$, we recall the Yoshino tensor product of matrix factorizations and find three of its variants.
Our theorem A is also stated and proved here.
In section $3$, properties of the various tensor products defined in the previous section are discussed. In section $4$, the class of summand-reducible polynomials is defined and the standard algorithm for matrix factorization of polynomials is improved on this class. Our theorem B is also stated and proved here. An example is provided to illustrate this result. The last section  is devoted to further research directions.
\section{Tensor products of matrix factorizations and their functoriality}
In this section, we first recall the definition of Yoshino's tensor product of matrix factorization denoted $\widehat{\otimes}$. Then, we find three bifunctorial variants of this product denoted by $\widehat{\otimes}'$, $\widehat{\otimes}''$ and $\widehat{\otimes}'''$. Next, we use the usual (standard) operations of tensor product and direct sum of matrices to construct a new product of matrix factorizations that we call \textit{multiplicative tensor product} of matrix factorizations denoted $\widetilde{\otimes}$. Moreover, we define its variant and then give some examples.\\
Under this section, unless otherwise stated, $R=K[[x]]$ and $S=K[[y]]$ where $x=x_{1},...,x_{r}$ and $y=y_{1},...,y_{s}$.

\subsection{Yoshino's tensor product of matrix factorizations and its variants} \label{subsec: Yoshino tens prod and variant}
Recall the following:
\begin{definition} \cite{yoshino1998tensor}  \label{defn matrix facto of polyn}   \\
An $n\times n$ \textbf{matrix factorization} of a power series $f\in \;R$ is a pair of $n$ $\times$ $n$ matrices $(\phi,\psi)$ such that
$\phi\psi=\psi\phi=fI_{n}$, where $I_{n}$ is the $n \times n$ identity matrix and the coefficients of $\phi$ and of $\psi$ are taken from $R$.
\end{definition}
Also recall ($\S 1$ of \cite{yoshino1998tensor}) the definition of
the category of matrix factorizations of a power series $f\in R=K[[x]]:=K[[x_{1},\cdots,x_{n}]]$ denoted by $MF(R,f)$ or $MF_{R}(f)$, (or even $MF(f)$ when there is no risk of confusion):\\
$\bullet$ The objects are the matrix factorizations of $f$.\\
$\bullet$ Given two matrix factorizations of $f$; $(\phi_{1},\psi_{1})$ and $(\phi_{2},\psi_{2})$ respectively of sizes $n_{1}$ and $n_{2}$, a morphism from $(\phi_{1},\psi_{1})$ to $(\phi_{2},\psi_{2})$ is a pair of matrices $(\alpha,\beta)$ each of size $n_{2}\times n_{1}$ which makes the following diagram commute \cite{yoshino1998tensor}:
$$\xymatrix@ R=0.6in @ C=.75in{K[[x]]^{n_{1}} \ar[r]^{\psi_{1}} \ar[d]_{\alpha} &
K[[x]]^{n_{1}} \ar[d]^{\beta} \ar[r]^{\phi_{1}} & K[[x]]^{n_{1}}\ar[d]^{\alpha}\\
K[[x]]^{n_{2}} \ar[r]^{\psi_{2}} & K[[x]]^{n_{2}}\ar[r]^{\phi_{2}} & K[[x]]^{n_{2}}}$$
That is,
$$\begin{cases}
 \alpha\phi_{1}=\phi_{2}\beta  \\
 \psi_{2}\alpha= \beta\psi_{1}
\end{cases}$$
More details on this category are found in chapter 2 of \cite{fomatati2019multiplicative}.
\\
\begin{definition} \cite{yoshino1998tensor} \label{defn Yoshino tensor prodt}
Let $X=(\phi,\psi)$ be an $n\times n$ matrix factorization of $f\in R$  and $X'=(\phi',\psi')$ an $m\times m$ matrix factorization of $g\in S$. These matrices can be considered as matrices over $L=K[[x,y]]$ and the \textbf{tensor product} $X\widehat{\otimes} X'$ is given by\\
(\(
\begin{bmatrix}
    \phi\otimes 1_{m}  &  1_{n}\otimes \phi'      \\
   -1_{n}\otimes \psi'  &  \psi\otimes 1_{m}
\end{bmatrix}
,
\begin{bmatrix}
    \psi\otimes 1_{m}  &  -1_{n}\otimes \phi'      \\
    1_{n}\otimes \psi'  &  \phi\otimes 1_{m}
\end{bmatrix}
\))\\
where each component is an endomorphism on $L^{n}\otimes L^{m}$.
\end{definition}
$X\widehat{\otimes} X'$ is an object of $MF_{L}(f+g)$ of size $2nm$ as proved in Lemma 2.1 of \cite{fomatati2019multiplicative}.
\begin{remark}
When $n=1$, we get a $1$ $\times$ $1$ matrix factorization of $f$, i.e., $f=[g][h]$ which is simply a factorization of $f$ in the classical sense. But in case $f$ is not reducible, this is not interesting, that's why we will mostly consider $n > 1$.
\end{remark}
\textbf{Variants of Yoshino's tensor product of matrix factorizations}\\
We now observe that in definition \ref{defn Yoshino tensor prodt}, each time we rotate the matrix on the left anticlockwise and the matrix on the right clockwise, we obtain a variant of the Yoshino tensor product.
\begin{definition} \label{defn Yoshino tensor prodt}
Let $X=(\phi,\psi)$ be an $n\times n$ matrix factorization of $f\in R$  and $X'=(\phi',\psi')$ an $m\times m$ matrix factorization of $g\in S$. These matrices can be considered as matrices over $L=K[[x,y]]$ and the \textbf{tensor products} $X\widehat{\otimes}' X'$, $X\widehat{\otimes}'' X'$ and $X\widehat{\otimes}''' X'$ are respectively given by\\
(\(
\begin{bmatrix}
  1_{n}\otimes \phi'  &  \psi\otimes 1_{m}  \\
    \phi\otimes 1_{m} & -1_{n}\otimes \psi'
\end{bmatrix}
,
\begin{bmatrix}
   1_{n}\otimes \psi' & \psi\otimes 1_{m}         \\
    \phi\otimes 1_{m} &  -1_{n}\otimes \phi'
\end{bmatrix}
\)),\\
(\(
\begin{bmatrix}
  \psi\otimes 1_{m} & -1_{n}\otimes \psi' \\
  1_{n}\otimes \phi' & \phi\otimes 1_{m}
\end{bmatrix}
,
\begin{bmatrix}
 \phi\otimes 1_{m} & 1_{n}\otimes \psi' \\
 -1_{n}\otimes \phi' & \psi\otimes 1_{m}
\end{bmatrix}
\)) and \\
(\(
\begin{bmatrix}
-1_{n}\otimes \psi' & \phi\otimes 1_{m}\\
\psi\otimes 1_{m}  & 1_{n}\otimes \phi'
\end{bmatrix}
,
\begin{bmatrix}
-1_{n}\otimes \phi' & \phi\otimes 1_{m} \\
 \psi\otimes 1_{m} & 1_{n}\otimes \psi'
\end{bmatrix}
\))\\
where each component is an endomorphism on $L^{n}\otimes L^{m}$.
\end{definition}
$X\widehat{\otimes}' X'$, $X\widehat{\otimes}'' X'$ and $X\widehat{\otimes}''' X'$ are objects of $MF_{L}(f+g)$, each of size $2nm$. The proof of this claim is completely similar to the proof of Lemma 2.1 of \cite{fomatati2019multiplicative}. Moreover, these three objects are mutually distinct as their definitions show.
\begin{definition}
  $\widehat{\otimes}'$, $\widehat{\otimes}''$ and $\widehat{\otimes}'''$ are respectively called the first, second and third variant of the Yoshino tensor product.
\end{definition}

\begin{proposition}
  $\widehat{\otimes}'$, $\widehat{\otimes}''$ and $\widehat{\otimes}'''$ are functorial operations.
\end{proposition}
\begin{proof}
  We prove that $\widehat{\otimes}'$ is a functorial operation. The proofs for $\widehat{\otimes}''$ and $\widehat{\otimes}'''$ are similar and therefore omitted.\\
  \textbf{Setting the stage:}
Let $X_{f}=(\phi,\psi)$, $X'_{f}=(\phi',\psi')$ and $X_{f}"=(\phi",\psi")$ be objects of $MF(K[[x]],f)$ respectively of sizes $n, n'$ and $n"$. Let $X_{g}=(\sigma,\rho)$, $X_{g}'=(\sigma',\rho')$ and $X_{g}"=(\sigma",\rho")$ be objects of $MF(K[[y]],g)$ respectively of sizes $m, m'$ and $m"$.

\begin{definition}
For morphisms $\zeta_{f}=(\alpha_{f}, \beta_{f}): X_{f}=(\phi,\psi) \rightarrow X_{f}'=(\phi',\psi')$  and $\zeta_{g}=(\alpha_{g}, \beta_{g}): X_{g}=(\sigma,\rho) \rightarrow X_{g}'=(\sigma',\rho')$ respectively in $MF(K[[x]],f)$ and $MF(K[[y]],g)$, we define $\zeta_{f}\widehat{\otimes}' X_{g}"$ by:
$$(
\begin{bmatrix}
  \alpha_{f}\otimes 1_{m"}  &          0      \\
    0               &   \alpha_{f}\otimes 1_{m"}
\end{bmatrix},
\begin{bmatrix}
     \beta_{f}\otimes 1_{m"}  &    0      \\
    0                 &   \beta_{f}\otimes 1_{m"}
\end{bmatrix}
)$$
We also define $X_{f}" \widehat{\otimes}' \zeta_{g}$
by
$$(
\begin{bmatrix}
  1_{n"}\otimes \alpha_{g}  &          0      \\
    0               &   1_{n"}\otimes \alpha_{g}
\end{bmatrix},
\begin{bmatrix}
     1_{n"}\otimes \beta_{g}  &    0      \\
    0                 &   1_{n"}\otimes \beta_{g}
\end{bmatrix}
)$$

\end{definition}

\begin{lemma}
\begin{enumerate}
\item
$\zeta_{f}\widehat{\otimes}' X_{g}"$: $X_{f}\widehat{\otimes}' X_{g}" =(\phi,\psi)\widehat{\otimes}' (\sigma",\rho")\rightarrow X_{f}'\widehat{\otimes}' X_{g}" =(\phi',\psi')\widehat{\otimes}' (\sigma",\rho")$ is a morphism in $MF(K[[x,y]],fg)$.
\item $X_{f}"\widehat{\otimes}' \zeta_{g}$: $X_{f}"\widehat{\otimes}' X_{g} =(\phi",\psi")\widehat{\otimes}' (\sigma,\rho)\rightarrow X_{f}"\widehat{\otimes}' X_{g}' =(\phi",\psi")\widehat{\otimes}' (\sigma',\rho')$ is a morphism in $MF(K[[x,y]],fg)$.
\end{enumerate}
\end{lemma}
\begin{proof}
  We omit the proof of this lemma because it is analogous to the proof we will give in the next section for the functoriality of the multiplicative tensor product of matrix factorizations.
\end{proof}
With all the above data, it is not difficult to see that for fixed $X_{f}$ and $X_{g}$, the tensor product $X_{f}\widehat{\otimes}' (-)$ (resp. $(-) \widehat{\otimes}' X_{g}$) defines the functor $MF_{K[[y]]}(g) \rightarrow MF_{K[[x,y]]}(fg)$ (resp. $MF_{K[[x]]}(f) \rightarrow MF_{K[[x,y]]}(fg)$).
\end{proof}
\subsection{Multiplicative tensor product of matrix factorization and its variant} \label{subsec: mult tens prodt n its variant}
In this subsection, all the properties we will prove for the multiplicative tensor product of matrix factorizations also hold for its variant and are proved similarly. So, we will omit the proofs involving the variant.\\
First recall that if $A$ (resp. $B$) is an $m\times n$ (resp. $p\times q$) matrix, then their direct sum $A\oplus B = \begin{bmatrix}
    A     &       0      \\
    0     &       B
\end{bmatrix}$, where the $0$ in the first line is a $m\times q$ matrix and the one in the second line is an $p\times n$ matrix. \\
Finally, recall that if $A$ (resp. $B$) is an $m\times n$ (resp. $p\times q$) matrix, then their tensor product $A\otimes B$ is the matrix obtained by replacing each entry $a_{ij}$ of $A$ with the matrix $a_{ij}B$. So, $A\otimes B$ is a $mp\times nq$ matrix.

\begin{definition} \label{defn of the multiplicative tensor product}
Let $X=(\phi,\psi)$ be a matrix factorization of $f\in K[[x]]$ of size $n$ and let $X'=(\phi',\psi')$ be a matrix factorization of $g\in K[[y]]$ of size $m$. Thus, $\phi,\psi,\phi' \,and \,\psi'$ can be considered as matrices over $L=K[[x,y]]$ and the \textbf{multiplicative tensor product} $X\widetilde{\otimes} X'$ is given by \\\\
\[((\phi\otimes\phi')\oplus (\phi\otimes\phi'), (\psi\otimes\psi')\oplus (\psi\otimes\psi'))=(
\begin{bmatrix}
    \phi\otimes\phi'  &          0      \\
    0                  &\phi\otimes\phi'
\end{bmatrix},
\begin{bmatrix}
    \psi\otimes\psi'  &    0      \\
    0                 &  \psi\otimes\psi'
\end{bmatrix}
)\]
\\\\
where each component is an endomorphism on $L^{n}\otimes_{L} L^{m}$.
\end{definition}
\begin{definition} Variant of $\widetilde{\otimes}$. \label{defn of the variant of the mult tens prod}\\
Let $X=(\phi,\psi)$ be a matrix factorization of $f\in K[[x]]$ of size $n$ and let $X'=(\phi',\psi')$ be a matrix factorization of $g\in K[[y]]$ of size $m$. Thus, $\phi,\psi,\phi' \,and \,\psi'$ can be considered as matrices over $L=K[[x,y]]$ and the \textbf{variant of the multiplicative tensor product} $X\widetilde{\otimes}' X'$ is given by \\\\
\[=(
\begin{bmatrix}
    0              &        \phi\otimes\phi'      \\
   \phi\otimes\phi' &          0
\end{bmatrix},
\begin{bmatrix}
       0          &  \psi\otimes\psi'        \\
\psi\otimes\psi'  &   0
\end{bmatrix}
)\]
\\\\
where each component is an endomorphism on $L^{n}\otimes_{L} L^{m}$.

\end{definition}

Observe that $X\widetilde{\otimes} X'$ and $X\widetilde{\otimes}' X'$ are distinct objects.

\begin{fact} \label{fact on tensor prod of aI and bI}
  Let $a$ and $b$ be two elements of the ring $K[[x_{1},\cdots,x_{n}]]$. \\
  Then, $aI_{n}\otimes bI_{m}=ab(I_{n}\otimes I_{m})$ where $\otimes$ is the tensor product of matrices.
\end{fact}
\begin{lemma} \label{lemma size of X tensor Y}
Let $X=(\phi,\psi)$ be an $n\times n$ matrix factorization of $f\in K[[x]]$ and let $X'=(\phi',\psi')$ be an $m\times m$ matrix factorization of $g\in K[[y]]$. Then,
$X\widetilde{\otimes} X'$ and $X\widetilde{\otimes}' X'$ are objects of $MF(K[[x,y]], fg)$ of size $2nm$.
\end{lemma}
\begin{proof}
We do it for $\widetilde{\otimes}$ and omit the proof for $\widetilde{\otimes}'$ because it is similar.
We have:\\\\
\[
\begin{bmatrix}
    \phi\otimes\phi'  &          0      \\
    0                  &\phi\otimes\phi'
\end{bmatrix}
\begin{bmatrix}
    \psi\otimes\psi'  &    0      \\
    0                 &  \psi\otimes\psi'
\end{bmatrix}\]
\\
\[=\begin{bmatrix}
    \phi\psi\otimes\phi'\psi'  &    0      \\
    0                 &  \phi\psi\otimes\phi'\psi'
\end{bmatrix}\]\\
\[=\begin{bmatrix}
    f1_{n}\otimes g1_{m}  &    0      \\
    0                 &  f1_{n}\otimes g1_{m}
\end{bmatrix},\,\,\,since\,\phi\psi=f1_{n}\,\,and\,\phi'\psi'=g1_{n}\]\\
\[=fg\begin{bmatrix}
    1_{n}\otimes 1_{m}  &    0      \\
    0                 &  1_{n}\otimes 1_{m}
\end{bmatrix}, by \,fact \,\ref{fact on tensor prod of aI and bI}.\]
\[=fg\cdot 1_{2nm}
\]
So, $X\widetilde{\otimes} X'$ is an object of $MF(fg)$ of size $2nm$ as claimed.
\end{proof}
\begin{example}
Consider the following $1\times 1$ matrix factorizations of $f=x^{3}$, $g=y^{5}$ and $h=z^{7}$:\\
$X=(x,x^{2})\in MF(x^{3})$, $X'=(y^{2},y^{3})\in MF(y^{5})$, $X''=(z^{3},z^{4})\in MF(z^{7})$.
We compute $X\widetilde{\otimes}X'$ and $(X\widetilde{\otimes}X')\widetilde{\otimes}X''$.\\\\
\[Y=X\widetilde{\otimes}X'=
(\begin{bmatrix}
    xy^{2}  &          0      \\
    0               &  xy^{2}
\end{bmatrix},
\begin{bmatrix}
    x^{2}y^{3}  &    0      \\
    0                 &  x^{2}y^{3}
\end{bmatrix})=(\phi_{Y}, \psi_{Y})\]\\
By lemma \ref{lemma size of X tensor Y}, $Y$ is a matrix factorization of $fg=x^{3}y^{5}$ of size $2(1)(1)=2$.
\\
Next, let $X''=(\phi_{X''},\psi_{X''})=(z^{3},z^{4})$, we compute:
\[Y\widetilde{\otimes} X''=(
\begin{bmatrix}
    \phi_{Y}\otimes \phi_{X''}     &          0      \\
                  0               &  \phi_{Y}\otimes \phi_{X''}
\end{bmatrix},
\begin{bmatrix}
   \psi_{Y}\otimes \psi_{X''}  &    0      \\
    0                 &  \psi_{Y}\otimes \psi_{X''}
\end{bmatrix}
)\]
\[=(
\begin{bmatrix}
    \begin{bmatrix}
    xy^{2}  &          0      \\
    0               &  xy^{2}
\end{bmatrix}\otimes z^{3}     &          0      \\
                  0               &  \begin{bmatrix}
    xy^{2}  &    0      \\
    0                 &  xy^{2}
\end{bmatrix}\otimes z^{3}
\end{bmatrix},
\begin{bmatrix}
   \begin{bmatrix}
    x^{2}y^{3}  &    0      \\
    0                 &  x^{2}y^{3}
\end{bmatrix}\otimes z^{4}  &    0      \\
    0                 &  \begin{bmatrix}
    x^{2}y^{3}  &    0      \\
    0                 &  x^{2}y^{3}
\end{bmatrix}\otimes z^{4}
\end{bmatrix}
)\]
\[=(\begin{bmatrix}
    xy^{2}z^{3}  &      0         &          0        &    0 \\
    0            &  xy^{2}z^{3}    &          0       &    0\\
    0            &     0          &      xy^{2}z^{3}  &   0\\
    0            &     0           &          0      & xy^{2}z^{3}
\end{bmatrix},\begin{bmatrix}
    x^{2}y^{3}z^{4}  &      0         &          0        &    0 \\
    0            &  x^{2}y^{3}z^{4}   &          0       &    0\\
    0            &     0          &      x^{2}y^{3}z^{4}  &   0\\
    0            &     0           &          0      & x^{2}y^{3}z^{4}
\end{bmatrix})\]


By lemma \ref{lemma size of X tensor Y}, $Y\widetilde{\otimes} X''\in MF(x^{3}y^{5}z^{7})=MF((fg)h)$ and $Y\widetilde{\otimes} X''$ is of size $2(1)(2)=4$ as $X''$ is of size $1$ and $Y$ is of size $2$.

\end{example}
\begin{example} We now give an example with $2\times 2$ matrix factorizations.\\
A straightforward computation shows that a $2 \times 2$ matrix factorization of $f=x^{2}+1$ is:\\
$
X=(\begin{bmatrix}
    x  &  -1      \\
    1  &  x
\end{bmatrix},
\begin{bmatrix}
    x  &  1      \\
    -1  & x
\end{bmatrix})=(\phi_{X},\psi_{X})
$\\\\
Hence, a matrix factorization of $g=-f=-(x^{2}+1)$ is:\\

$
Y=(\begin{bmatrix}
    -x  &  1      \\
    -1  &  -x
\end{bmatrix},
\begin{bmatrix}
    x  &  1      \\
    -1  & x
\end{bmatrix})=(\phi_{Y},\psi_{Y})
$\\\\
\[ X\widetilde{\otimes} Y=(\begin{bmatrix}
    \phi_{X}\otimes \phi_{Y}  &  0      \\
    0  &  \phi_{X}\otimes \phi_{Y}
\end{bmatrix},
\begin{bmatrix}
    \psi_{X}\otimes \psi_{Y}  &  0      \\
    0  & \psi_{X}\otimes \psi_{Y}
\end{bmatrix})
\]\\
\begin{gather*}
  \setlength{\arraycolsep}{1.0\arraycolsep}
  \text{\footnotesize$\displaystyle
  i.e.,\, X\widetilde{\otimes} Y =(\begin{bmatrix}
    \begin{bmatrix}
    -x^{2} & x & x  &  -1      \\
    -x  &  -x^{2} & 1 & x \\
    -x  & 1  &   -x^{2}& x \\
    -1  &-x  &   -x   &  -x^{2}
\end{bmatrix}  &  0      \\
    0  &  \begin{bmatrix}
    -x^{2} & x & x  &  -1      \\
    -x  &  -x^{2} & 1 & x \\
    -x  & 1  &   -x^{2}& x \\
    -1  &-x  &   -x   &  -x^{2}
\end{bmatrix}
\end{bmatrix},
\begin{bmatrix}
    \begin{bmatrix}
    x^{2} & x & x  &   1      \\
    -x  &  x^{2} & -1 &  x \\
    -x  & -1    &  x^{2}& x \\
    1  &-x  &   -x   &  x^{2}
\end{bmatrix}  &  0      \\
    0  &  \begin{bmatrix}
    x^{2} & x & x  &   1      \\
    -x  &  x^{2} & -1 &  x \\
    -x  & -1    &  x^{2}& x \\
    1  &-x  &   -x   &  x^{2}
\end{bmatrix}
\end{bmatrix})$}
\end{gather*}
\\
By lemma \ref{lemma size of X tensor Y}, $X\widetilde{\otimes} Y \in MF(fg)$ and $X\widetilde{\otimes }Y$ is of size $2(2)(2)=8$ as $X$ and $Y$ are both of size $2$.
\end{example}

%

%



\subsection{Functoriality of the operations $\widetilde{\otimes}$ and $\widetilde{\otimes}'$ }
This subsection is entirely devoted to the discussion of the bifunctoriality of $\widetilde{\otimes}$. The proof of the bifunctoriality of $\widetilde{\otimes}'$ is omitted because it is similar to that of $\widetilde{\otimes}$. \\
Let $X_{f}$, $X_{f}'$, $X_{f}"$, $X_{g}$, $X_{g}'$ and $X_{g}"$ be defined as in subsection \ref{subsec: Yoshino tens prod and variant}.

\begin{definition}\label{defn zeta is a bifuntor}
For morphisms $\zeta_{f}=(\alpha_{f}, \beta_{f}): X_{f}=(\phi,\psi) \rightarrow X_{f}'=(\phi',\psi')$  and $\zeta_{g}=(\alpha_{g}, \beta_{g}): X_{g}=(\sigma,\rho) \rightarrow X_{g}'=(\sigma',\rho')$ respectively in $MF(K[[x]],f)$ and $MF(K[[y]],g)$,\\ we define $\zeta_{f}\widetilde{\otimes} \zeta_{g}: X_{f}\widetilde{\otimes}X_{g} =(\phi,\psi)\widetilde{\otimes} (\sigma,\rho)\rightarrow X_{f}'\widetilde{\otimes}X_{g}' =(\phi',\psi')\widetilde{\otimes} (\sigma',\rho')$
by
$$(
\begin{bmatrix}
  \alpha_{f}\otimes \alpha_{g}  &          0      \\
    0               &   \alpha_{f}\otimes \alpha_{g}
\end{bmatrix},
\begin{bmatrix}
     \beta_{f}\otimes \beta_{g}  &    0      \\
    0                 &   \beta_{f}\otimes \beta_{g}
\end{bmatrix}
).$$
$\zeta_{f}\widetilde{\otimes}' \zeta_{g}$ is defined the same way.
\end{definition}

\begin{lemma}\label{zeta is a morphism in both arguments}
\begin{enumerate}
  \item $\zeta_{f}\widetilde{\otimes} \zeta_{g}$: $X_{f}\widetilde{\otimes}X_{g} =(\phi,\psi)\widetilde{\otimes} (\sigma,\rho)\rightarrow X_{f}'\widetilde{\otimes}X_{g}' =(\phi',\psi')\widetilde{\otimes} (\sigma',\rho')$ is a morphism in $MF(K[[x,y]],fg)$.
  \item $\zeta_{f}\widetilde{\otimes}' \zeta_{g}$: $X_{f}\widetilde{\otimes}'X_{g} =(\phi,\psi)\widetilde{\otimes}' (\sigma,\rho)\rightarrow X_{f}'\widetilde{\otimes}'X_{g}' =(\phi',\psi')\widetilde{\otimes}' (\sigma',\rho')$ is a morphism in $MF(K[[x,y]],fg)$.
\end{enumerate}
\end{lemma}
\begin{proof}
We do it only for $\widetilde{\otimes}$ since the proof for $\widetilde{\otimes}'$ is very similar. We need to show that the following diagram commutes:

 $$\xymatrix@ R=0.8in @ C=1.10in{K[[x,y]]^{2nm} \ar[r]^{\begin{bmatrix}
    \psi\otimes\rho  &    0      \\
    0                 &  \psi\otimes\rho
\end{bmatrix}} \ar[d]_{\begin{bmatrix}
  \alpha_{f} \otimes \alpha_{g}  &          0      \\
    0               &  \alpha_{f}\otimes \alpha_{g}
\end{bmatrix}} &
K[[x,y]]^{2nm} \ar[d]^{\begin{bmatrix}
    \beta_{f} \otimes \beta_{g}   &    0      \\
    0                 &    \beta_{f} \otimes \beta_{g}
\end{bmatrix}} \ar[r]^{\begin{bmatrix}
    \phi\otimes\sigma  &          0      \\
    0                  &\phi\otimes\sigma
\end{bmatrix}} & K[[x,y]]^{2nm}\ar[d]^{\begin{bmatrix}
  \alpha_{f} \otimes \alpha_{g}  &          0      \\
    0               &  \alpha_{f}\otimes \alpha_{g}
\end{bmatrix}}\\
K[[x,y]]^{2n'm'} \ar[r]_{\begin{bmatrix}
    \psi'\otimes\rho'  &    0      \\
    0                 &  \psi'\otimes\rho'
\end{bmatrix}} & K[[x,y]]^{2n'm'}\ar[r]_{\begin{bmatrix}
    \phi'\otimes\sigma'  &          0      \\
    0                  &\phi'\otimes\sigma'
\end{bmatrix}} & K[[x,y]]^{2n'm'}}$$

 viz. both the left and the right squares in the foregoing diagram commute.\\
 $\bullet$ The commutativity of the right square and the left square are respectively expressed by the following equalities:
 \[\begin{bmatrix}
  \alpha_{f} \otimes \alpha_{g}  &          0      \\
    0               &  \alpha_{f}\otimes \alpha_{g}
\end{bmatrix}\begin{bmatrix}
    \phi\otimes\sigma  &          0      \\
    0                  &\phi\otimes\sigma
\end{bmatrix}=\begin{bmatrix}
    \phi'\otimes\sigma'  &          0      \\
    0                  &\phi'\otimes\sigma'
\end{bmatrix}\begin{bmatrix}
    \beta_{f} \otimes \beta_{g}   &    0      \\
    0                 &    \beta_{f} \otimes \beta_{g}
\end{bmatrix} \]
and\\
 \[\begin{bmatrix}
    \beta_{f} \otimes \beta_{g}   &    0      \\
    0                 &    \beta_{f} \otimes \beta_{g}
\end{bmatrix}\begin{bmatrix}
    \psi\otimes\rho  &    0      \\
    0                 &  \psi\otimes\rho
\end{bmatrix}=\begin{bmatrix}
    \psi'\otimes\rho'  &    0      \\
    0                 &  \psi'\otimes\rho'
\end{bmatrix}\begin{bmatrix}
  \alpha_{f} \otimes \alpha_{g}  &          0      \\
    0               &  \alpha_{f}\otimes \alpha_{g}
\end{bmatrix} \]


%
i.e., all we need to show is the pair of equalities:
\\
$$\begin{cases}
\alpha_{f}\phi\otimes \alpha_{g}\sigma=\phi'\beta_{f}\otimes\sigma'\beta_{g} \cdots (1)\\
\beta_{f}\psi \otimes \beta_{g}\rho= \psi'\alpha_{f}\otimes\rho'\alpha_{g} \cdots (2)
\end{cases}$$

Now by hypothesis, $\zeta_{f}=(\alpha_{f},\beta_{f}): X_{f}=(\phi,\psi) \rightarrow X_{f}'=(\phi',\psi')$ and $\zeta_{g}=(\alpha_{g},\beta_{g}): X_{g}=(\sigma,\rho) \rightarrow X_{g}'=(\sigma',\rho')$ are morphisms, meaning that the following diagrams commute
$$\xymatrix@ R=0.6in @ C=.75in{K[[x]]^{n} \ar[r]^{\psi} \ar[d]_{\alpha_{f}} &
K[[x]]^{n} \ar[d]^{\beta_{f}} \ar[r]^{\phi} & K[[x]]^{n}\ar[d]^{\alpha_{f}}\\
K[[x]]^{n'} \ar[r]^{\psi'} & K[[x]]^{n'}\ar[r]^{\phi'} & K[[x]]^{n'}}$$

and \\
$$\xymatrix@ R=0.6in @ C=.75in{K[[y]]^{m} \ar[r]^{\rho} \ar[d]_{\alpha_{g}} &
K[[y]]^{m} \ar[d]^{\beta_{g}} \ar[r]^{\sigma} & K[[y]]^{m}\ar[d]^{\alpha_{g}}\\
K[[y]]^{m'} \ar[r]^{\rho'} & K[[y]]^{m'}\ar[r]^{\sigma'} & K[[y]]^{m'}}$$

  That is,
$$\begin{cases}
 \alpha_{f}\phi=\phi'\beta_{f} \cdots (i) \\
 \psi'\alpha_{f}= \beta_{f}\psi \cdots (ii)
\end{cases}$$
and \\
$$\begin{cases}
 \alpha_{g}\sigma=\sigma'\beta_{g} \cdots (i') \\
 \rho'\alpha_{g}= \beta_{g}\rho \cdots (ii')
\end{cases}$$
Now considering $(i)$ and $(i')$, we immediately see that equality $(1)$ holds. Similarly, $(ii)$ and $(ii')$ yield $(2)$. \\

So, $\zeta_{f} \widetilde{\otimes} \zeta_{g}$ is a morphism in $MF(K[[x,y]], fg)$.

\end{proof}
We can now state the following result.
\begin{theorem} \label{zeta is a bifunctor}
\begin{enumerate}
\item If $X$ (resp. $Y$) is a matrix factorization of $f$ (resp. $g$). Then, there is a tensor product $\widetilde{\otimes}$ of matrix factorizations which produces a matrix factorization $X \widetilde{\otimes} Y$ of the product $fg\in K[[x_{1},...,x_{r},y_{1},...,y_{s}]]$. $\widetilde{\otimes}$ is called the multiplicative tensor product of matrix factorizations. Moreover, $\widetilde{\otimes}$ has a functorial variant $\widetilde{\otimes}'$ having the same effect on $X$ and $Y$.

\item The multiplicative tensor product $(-) \widetilde{\otimes} (-):MF(K[[x]],f)\times MF(K[[y]],g)\rightarrow MF(K[[x,y]],fg)$ is a bifunctor. Furthermore, its variant $\widetilde{\otimes}'$ is also a bifunctor.
\end{enumerate}
\end{theorem}
\begin{proof}
\begin{enumerate}
  \item This is exactly what we proved above in subsection \ref{subsec: mult tens prodt n its variant}.
  \item
 We focus only on the proof for $\widetilde{\otimes}$ and omit that for $\widetilde{\otimes}'$ because the proofs are similar.\\
In order to ease our computations, let's write $F=(-)\widetilde{\otimes} (-)$. We show that $F$ is a bifunctor.

\newpage We have:

$\;\; (-)\widetilde{\otimes} (-):\;\;\;\;\;\;\;\;\;\;\;\;\; MF(f)\times MF(g)\,\,\,\,\,\,\,\,\,\,\longrightarrow \,\,\,\,\,\,\,\,\,\,\,\,\,\,\,\,MF(fg)$
$$\xymatrix @ R=0.4in @ C=.33in
{&(X_{f}\;\;\;\;\;\;{,} \ar[d]_{\zeta_{f}}\ar @{}[dr]& X_{g}) \ar[rrr]^-{} \ar[d]_{\zeta_{g}}
\ar @{}[dr] &&& X_{f}\widetilde{\otimes} X_{g} \ar[d]^{\zeta_{f}\widetilde{\otimes} \zeta_{g}:= (\alpha,\beta)}\\
&(X_{f}'\;\;\;\;\;\; {,}& X_{g}') &\ar[r]&& X_{f}'\widetilde{\otimes} X_{g}'}$$
$$\xymatrix @ R=0.4in @ C=.3in
{&\ar[d]_{\zeta_{f}'}\ar @{}[dr]& \ar[d]_{\zeta_{g}'}
\ar @{}[dr] &&&\ar[d]^{\zeta_{f}'\widetilde{\otimes} \zeta_{g}':= (\alpha',\beta')}\\
&(X_{f}" \;\;\;\;\;\;{,}& X_{g}") &\ar[r]&& X_{f}"\widetilde{\otimes} X_{g}"}$$

We showed in lemma \ref{zeta is a morphism in both arguments} that $\zeta_{f}\widetilde{\otimes} \zeta_{g}:= (\alpha,\beta)$ is a morphism in
$MF(K[[x,y]],fg)$, where $$(\alpha,\beta)=(
\begin{bmatrix}
  \alpha_{f}\otimes \alpha_{g}  &          0      \\
    0               &   \alpha_{f}\otimes \alpha_{g}
\end{bmatrix},
\begin{bmatrix}
     \beta_{f}\otimes \beta_{g}  &    0      \\
    0                 &   \beta_{f}\otimes \beta_{g}
\end{bmatrix}
).$$ Similarly, if $\zeta_{f}':=(\alpha_{f}',\beta_{f}')$ and $\zeta_{g}':=(\alpha_{g}',\beta_{g}')$ then $\zeta_{f}' \widetilde{\otimes} \zeta_{g}'=(\alpha',\beta')$ where $$(\alpha',\beta')=(
\begin{bmatrix}
  \alpha_{f}'\otimes \alpha_{g}'  &          0      \\
    0               &   \alpha_{f}'\otimes \alpha_{g}'
\end{bmatrix},
\begin{bmatrix}
     \beta_{f}'\otimes \beta_{g}'  &    0      \\
    0                 &   \beta_{f}'\otimes \beta_{g}'
\end{bmatrix}
).$$ It now remains to show the composition and the identity axioms.\\\\
\textit{Identity Axiom}: \\
We show that $F(id_{(X_{f},X_{g})})=id_{F(X_{f},X_{g})}$.\\
Now, $F(id_{(X_{f},X_{g})})=F(id_{X_{f}},id_{X_{g}}):=id_{X_{f}}\widetilde{\otimes} id_{X_{g}}: X_{f} \widetilde{\otimes}X_{g} \rightarrow X_{f} \widetilde{\otimes}X_{g}$.\\
And by definition \ref{defn zeta is a bifuntor}, $id_{X_{f}}\widetilde{\otimes} id_{X_{g}}$ is the pair of matrices \\
 $(
\begin{bmatrix}
  I_{n}\otimes I_{m}  &          0      \\
    0               &    I_{n}\otimes I_{m}
\end{bmatrix},
\begin{bmatrix}
 I_{n}\otimes I_{m}  &          0      \\
    0               &    I_{n}\otimes I_{m}
\end{bmatrix})\,\,\,\,\,\,\,\,\,\,\,\dag$ \\

Next, we compute $id_{F(X_{f},X_{g})}=id_{X_{f} \widetilde{\otimes}X_{g}}:X_{f} \widetilde{\otimes}X_{g} \rightarrow X_{f} \widetilde{\otimes}X_{g}$.\\
By definition of a morphism in the category $MF(fg)$, we know that \\

$id_{X_{f} \widetilde{\otimes}X_{g}}:=(
\begin{bmatrix}
  I_{nm}  &          0      \\
    0               &    I_{nm}
\end{bmatrix},
\begin{bmatrix}
 I_{nm}  &          0      \\
    0               &    I_{nm}
\end{bmatrix})\,\,\,\,\,\,\,\,\,\,\,\dag\dag
$\\

Since $I_{n}\otimes I_{m}= I_{nm}$, we see that $\dag$ and $\dag\dag$ are the same, therefore $F(id_{(X_{f},X_{g})})=id_{F(X_{f},X_{g})}$ as desired.\\\\
\textit{Composition Axiom}:\\
Consider the situation:
$$\xymatrix@R=4in{X_{f}\ar[r]^{\zeta_{f}} & X_{f}'\ar[r]^{\zeta_{f}'}& X_{f}"}$$ $$\xymatrix@R=4in{X_{g}\ar[r]^{\zeta_{g}} & X_{g}'\ar[r]^{\zeta_{g}'}& X_{g}"}$$
$$\xymatrix@R=4in{X_{f} \widetilde{\otimes}X_{g}\ar[r]^{F(\zeta_{f},\zeta_{g})} & X_{f}'\widetilde{\otimes}X_{g}'\ar[r]^{F(\zeta_{f}',\zeta_{g}')}& X_{f}"\widetilde{\otimes} X_{g}"}$$

We need to show $F(\zeta_{f}'\circ \zeta_{f},\zeta_{g}'\circ \zeta_{g})=F(\zeta_{f}',\zeta_{g}')\circ F(\zeta_{f},\zeta_{g})$. \\
Now, $\zeta_{f}'\circ \zeta_{f}= (\alpha_{f}'\alpha_{f},\beta_{f}'\beta_{f})$ and $\zeta_{g}'\circ \zeta_{g}=(\alpha_{g}'\alpha_{g},\beta_{g}'\beta_{g})$.\\\\
Thanks to definition \ref{defn zeta is a bifuntor}, we obtain:
$$(\zeta_{f}'\circ \zeta_{f})\widetilde{\otimes}(\zeta_{g}'\circ \zeta_{g})=(
\begin{bmatrix}
  \alpha_{f}'\alpha_{f}\otimes \alpha_{g}'\alpha_{g}  &          0      \\
    0               &   \alpha_{f}'\alpha_{f}\otimes \alpha_{g}'\alpha_{g}
\end{bmatrix},
\begin{bmatrix}
     \beta_{f}'\beta_{f}\otimes \beta_{g}'\beta_{g}  &    0      \\
    0                 &   \beta_{f}'\beta_{f}\otimes \beta_{g}'\beta_{g}
\end{bmatrix}
)\,\,\,\,\,\,\,\,\,\,\,\, \ddag'$$

Next, \\\\
$(\zeta_{f}'\widetilde{\otimes}\zeta_{g}')\circ (\zeta_{f}\widetilde{\otimes} \zeta_{g})\\
=(
\begin{bmatrix}
  \alpha_{f}'\otimes \alpha_{g}'  &          0      \\
    0               &   \alpha_{f}'\otimes \alpha_{g}'
\end{bmatrix},
\begin{bmatrix}
     \beta_{f}'\otimes \beta_{g}'    &    0      \\
    0                              &   \beta_{f}'\otimes \beta_{g}'
\end{bmatrix}
)\circ
(
\begin{bmatrix}
  \alpha_{f}\otimes \alpha_{g}  &          0      \\
    0               &   \alpha_{f}\otimes \alpha_{g}
\end{bmatrix},
\begin{bmatrix}
     \beta_{f}\otimes \beta_{g}    &    0      \\
    0                              &   \beta_{f}\otimes \beta_{g}
\end{bmatrix}
)\\
=(
\begin{bmatrix}
  \alpha_{f}'\alpha_{f}\otimes \alpha_{g}'\alpha_{g}  &          0      \\
    0               &   \alpha_{f}'\alpha_{f}\otimes \alpha_{g}'\alpha_{g}
\end{bmatrix},
\begin{bmatrix}
     \beta_{f}'\beta_{f}\otimes \beta_{g}'\beta_{g}  &    0      \\
    0                 &   \beta_{f}'\beta_{f}\otimes \beta_{g}'\beta_{g}
\end{bmatrix}
)\,\,\,\,\,\,\,\,\,\,\,\, \ddag\ddag'$

From $\ddag'$ and $\ddag\ddag'$, we see that $F(\zeta_{f}'\circ \zeta_{f},\zeta_{g}'\circ \zeta_{g})=F(\zeta_{f}',\zeta_{g}')\circ F(\zeta_{f},\zeta_{g})$.
Thus, $(-)\widetilde{\otimes}(-)$ is a bifunctor.

\end{enumerate}
\end{proof}

\section{Properties of the multiplicative tensor product of matrix factorizations and of its variant}
In this section, we prove that $\widetilde{\otimes}$ and its variant are associative, commutative and distributive.
\\
We denote by $X_{1}=(\phi_{1},\psi_{1})$ (resp. $X_{2}=(\phi_{2},\psi_{2})$) an $(n_{1}\times n_{1})$ (resp. $(n_{2}\times n_{2})$) matrix factorization of $f\in K[[x]]$. We also let $X'=(\phi',\psi')$ (resp. $X''=(\phi'',\psi'')$) denotes a $(p\times p)$ (resp. $(m\times m)$) matrix factorization of $g\in K[[y]]$ (resp. of $h\in K[[z]]:= K[[z_{1},\cdots,z_{l}]]$).
$X=(\phi,\psi)$ will also be an $r\times r$ matrix factorization of $f\in K[[x]]$.

\subsection{Associativity, commutativity and distributivity of $\widetilde{\otimes}$ and $\widetilde{\otimes}'$}

\begin{proposition} (Associativity) \label{Assoc of new tensor product}\\
There are identities
\begin{enumerate}
\item  $(X\widetilde{\otimes}X')\widetilde{\otimes}X''=X\widetilde{\otimes}(X'\widetilde{\otimes}X'')$ in $MF(fgh)$.
\item  $(X\widetilde{\otimes}'X')\widetilde{\otimes}'X''=X\widetilde{\otimes}'(X'\widetilde{\otimes}'X'')$ in $MF(fgh)$.
\end{enumerate}
\end{proposition}
\begin{proof}
The desired identities follow from the fact that the standard tensor product for matrices is associative.
\end{proof}
Before stating the next proposition, it is worthwhile recalling (cf. section 3.1 \cite{henderson1981vec}) that given two matrices $A$ and $B$, the tensor products $A \otimes B$ and $B \otimes A$ are \textbf{permutation equivalent}. That is, there exist permutation matrices $P$ and $Q$ (so called commutation matrices) such that:
$A \otimes B = P (B \otimes A) Q$. If $A$ and $B$ are square matrices, then $A \otimes B$ and $B\otimes A$ are even \textbf{permutation similar}, meaning we can take $P = Q^{T}$.\\
To be more precise \cite{henderson1981vec}, if $A$ is a $p\times q$ matrix and $B$ is an $r\times s$ matrix, then $$B \otimes A = S_{p,r} (A \otimes B) S_{q,s}^{T}$$
where, $$ S_{m,n}=\sum_{i=1}^{m}( e_{i}^{T}\otimes I_{n}\otimes e_{i}) = \sum_{j=1}^{n}( e_{j}\otimes I_{m}\otimes e_{j}^{T})$$
$I_{n}$ is the $n\times n$ identity matrix and $e_{i}$ is the $i^{th}$ unit vector. $ S_{m,n}$ is the \textbf{perfect shuffle} permutation matrix.\\
We use $2\times 2$ matrices to illustrate the fact that $A \otimes B$ and $B \otimes A$ are permutation similar.\\
Let \(A=
\begin{bmatrix}
    a  &  b      \\
    c  &  d
\end{bmatrix}\) and \(B=
\begin{bmatrix}
    e  &  f      \\
    g  &  h
\end{bmatrix}.\)
Then \(A\otimes B=
\begin{bmatrix}
     a  &  b      \\
    c  &  d
\end{bmatrix}
\otimes
\begin{bmatrix}
     e  &  f      \\
    g  &  h
\end{bmatrix} =\begin{bmatrix}
     ae  &  af  & be & bf    \\
     ag  &  ah  & bg & bh    \\
     ce  &  cf  & de & df    \\
     cg  &  ch  & dg & dh
\end{bmatrix}\) and
\(B\otimes A=
\begin{bmatrix}
   e  &  f      \\
    g  &  h
\end{bmatrix}
\otimes
\begin{bmatrix}
  a  &  b      \\
    c  &  d
\end{bmatrix} =\begin{bmatrix}
     ea  &  eb  & fa & fb    \\
     ec  &  ed  & fc & fd    \\
     ga  &  gb  & ha & hb    \\
     gc  &  gd  & hc & hd
\end{bmatrix}\) and we have: \\
\(B\otimes A= \begin{bmatrix}
     ea  &  eb  & fa & fb    \\
     ec  &  ed  & fc & fd    \\
     ga  &  gb  & ha & hb    \\
     gc  &  gd  & hc & hd
\end{bmatrix} \leftrightsquigarrow_{c_{2}\leftrightarrow c_{3}} \begin{bmatrix}
     ea  &  fa  & eb & fb    \\
     ec  &  fc  & ed & fd    \\
     ga  &  ha  & gb & hb    \\
     gc  &  hc  & gd & hd
\end{bmatrix}
 \leftrightsquigarrow_{r_{2}\leftrightarrow r_{3}} \begin{bmatrix}
     ea  &  fa  & eb & fb    \\
     ga  &  ha  & gb & hb    \\
     ec  &  fc  & ed & fd    \\
     gc  &  hc  & gd & hd
\end{bmatrix} = A\otimes B \)
The $\leftrightsquigarrow_{c_{2}\leftrightarrow c_{3}}$ (respectively $\leftrightsquigarrow_{r_{2}\leftrightarrow r_{3}}$) mean that the second and third column (respectively the second and third row) have been interchanged.
The commutativity of $\widetilde{\otimes}$ is up to isomorphism. This isomorphism comes from the permutation similarity\footnote{Recall that all the matrices involved in a matrix factorization are square matrices by definition, this justifies the fact that we talk of permutation similarity instead of permutation equivalence.} of the matrices involved.
\begin{proposition} (commutativity) \label{prop commutativity of new tensor product}\\
For matrix factorizations $X\in MF(f)$ and $X'\in MF(g)$, we have
\begin{enumerate}
  \item $X\widetilde{\otimes} X'\cong X' \widetilde{\otimes} X \,in\,MF(fg).$
  \item $X\widetilde{\otimes}' X'\cong X' \widetilde{\otimes}' X \,in\,MF(fg).$
\end{enumerate}

\end{proposition}
\begin{proof} Since the proof of the two claims are similar, we only prove the first one.\\
We first prove that there is a morphism from the matrix factorization $X\widetilde{\otimes} X'$ to the matrix factorization $X'\widetilde{\otimes} X$. We know that:\\
$X\widetilde{\otimes} X' =(\begin{bmatrix}
  \phi\otimes \phi' &          0      \\
    0               &   \phi\otimes \phi'
\end{bmatrix},\begin{bmatrix}
  \psi\otimes \psi' &          0      \\
    0               &   \psi\otimes \psi'
\end{bmatrix})$;
$X'\widetilde{\otimes} X =(\begin{bmatrix}
  \phi'\otimes \phi &          0      \\
    0               &   \phi'\otimes \phi
\end{bmatrix},\begin{bmatrix}
  \psi'\otimes \psi &          0      \\
    0               &   \psi'\otimes \psi
\end{bmatrix})$
 Recall that $X$ and $X'$ are respectively of sizes $r$ and $p$. By definition of a morphism in $MF(fg)$, we find a pair of matrices $(\delta,\beta)$ such that the following diagram commutes:
$$\xymatrix@ R=0.6in @ C=.95in{K[[x,y]]^{2rp} \ar[r]^{\begin{bmatrix}
  \psi\otimes \psi' &          0      \\
    0               &   \psi\otimes \psi'
\end{bmatrix}} \ar[d]_{\delta} &
K[[x,y]]^{2rp} \ar[d]^{\beta} \ar[r]^{\begin{bmatrix}
  \phi\otimes \phi' &          0      \\
    0               &   \phi\otimes \phi'
\end{bmatrix}} & K[[x,y]]^{2rp}\ar[d]^{\delta}\\
K[[x,y]]^{2rp} \ar[r]_{\begin{bmatrix}
  \psi'\otimes \psi &          0      \\
    0               &   \psi'\otimes \psi
\end{bmatrix}} & K[[x,y]]^{2rp}\ar[r]_{\begin{bmatrix}
  \phi'\otimes \phi &          0      \\
    0               &   \phi'\otimes \phi
\end{bmatrix}} & K[[x,y]]^{2rp}}$$
 It suffices to choose $(\delta=\begin{bmatrix}
  \phi'\otimes \phi &          0      \\
    0               &   \phi'\otimes \phi
\end{bmatrix},\beta=\begin{bmatrix}
  \phi\otimes \phi' &          0      \\
    0               &   \phi\otimes \phi'
\end{bmatrix})$, for the above diagram to commute. In fact, the commutativity of the right square is immediate. As for the left square, we need to verify the following equality:
 $$\begin{bmatrix}
  \psi'\otimes \psi &          0      \\
    0               &   \psi'\otimes \psi
\end{bmatrix}\begin{bmatrix}
  \phi'\otimes \phi &          0      \\
    0               &   \phi'\otimes \phi
\end{bmatrix}= \begin{bmatrix}
  \phi\otimes \phi' &          0      \\
    0               &   \phi\otimes \phi'
\end{bmatrix}\begin{bmatrix}
  \psi\otimes \psi' &          0      \\
    0               &   \psi\otimes \psi'
\end{bmatrix}$$

i.e.,\\
 $(\psi'\otimes \psi)(\phi'\otimes \phi)=(\phi\otimes \phi')(\psi\otimes \psi')$ viz. $\psi'\phi'\otimes \psi\phi=\phi\psi\otimes \phi'\psi'$, i.e., $g\cdot I_{p} \otimes f\cdot I_{r}=f\cdot I_{r}\otimes g\cdot I_{p}$, i.e., $gfI_{p} \otimes I_{r}=fgI_{r} \otimes I_{p} $ which is true since $fg=gf$ and $I_{r} \otimes I_{p}=I_{p} \otimes I_{r}$.\\
 So, there is a map from $X\widetilde{\otimes} X'$ to $X'\widetilde{\otimes} X$.\\
Secondly, we prove the isomorphism:\\
In fact, $$X\widetilde{\otimes} X' =(\begin{bmatrix}
  \phi\otimes \phi' &          0      \\
    0               &   \phi\otimes \phi'
\end{bmatrix},\begin{bmatrix}
  \psi\otimes \psi' &          0      \\
    0               &   \psi\otimes \psi'
\end{bmatrix})$$
$$\cong (\begin{bmatrix}
  \phi'\otimes \phi &          0      \\
    0               &   \phi'\otimes \phi
\end{bmatrix},\begin{bmatrix}
  \psi'\otimes \psi &          0      \\
    0               &   \psi'\otimes \psi
\end{bmatrix})=X'\widetilde{\otimes} X$$
The "$ \cong $" in this proof is due to the fact that $\phi\otimes \phi'$ (respectively $\psi\otimes \psi'$) and $\phi'\otimes \phi$ (respectively $\psi'\otimes \psi$) are permutation similar.
\end{proof}
\begin{proposition} (Distributivity)\\
If $X_{1}$ and $X_{2}$ are matrix factorizations (of $f\in K[[x]]$) of the same size, then there are natural isomorphisms
\begin{enumerate}
  \item $(X_{1}\oplus X_{2})\widetilde{\otimes} X'\cong(X_{1}\widetilde{\otimes} X')\oplus (X_{2}\widetilde{\otimes} X').$
  \item $ X' \widetilde{\otimes}(X_{1}\oplus X_{2})\cong(X'\widetilde{\otimes}X_{1})\oplus (X'\widetilde{\otimes}X_{2}).$
   \item $(X_{1}\oplus X_{2})\widetilde{\otimes}' X'\cong(X_{1}\widetilde{\otimes}' X')\oplus (X_{2}\widetilde{\otimes}' X').$
  \item $ X' \widetilde{\otimes}'(X_{1}\oplus X_{2})\cong(X'\widetilde{\otimes}'X_{1})\oplus (X'\widetilde{\otimes}'X_{2}).$
\end{enumerate}

\end{proposition}
\begin{proof} We prove only the first two claims, since the last two claims are proved similarly.
\begin{enumerate}
  \item
$$(X_{1}\widetilde{\otimes} X')\oplus (X_{2}\widetilde{\otimes} X')$$
\[=(\begin{bmatrix}
  \phi_{1}\otimes \phi' &          0      \\
    0               &   \phi_{1}\otimes \phi'
\end{bmatrix}, \begin{bmatrix}
  \psi_{1}\otimes \psi'  &          0      \\
    0               &   \psi_{1}\otimes \psi'
\end{bmatrix})\oplus (\begin{bmatrix}
  \phi_{2}\otimes \phi' &          0      \\
    0               &   \phi_{2}\otimes \phi'
\end{bmatrix}, \begin{bmatrix}
  \psi_{2}\otimes \psi'   &          0      \\
    0               &   \psi_{2}\otimes \psi'
\end{bmatrix})\]\\
\[=(\begin{bmatrix}
  \phi_{1}\otimes \phi' &   0  &  0   &    0      \\
  0   & \phi_{1}\otimes \phi' &  0   &    0  \\
   0   &  0   & \phi_{2}\otimes \phi' &    0 \\
    0    &  0  &  0     &   \phi_{2}\otimes \phi'
\end{bmatrix}, \begin{bmatrix}
  \psi_{1}\otimes \psi' &   0  &  0   &    0      \\
  0   & \psi_{1}\otimes \psi' &  0   &    0  \\
   0   &  0   & \psi_{2}\otimes \psi' &    0 \\
    0    &  0  &  0     &   \psi_{2}\otimes \psi'
\end{bmatrix})\,\,\,\,\,\,\,\,\cdots (\sharp)\]\\

Let $P=\begin{bmatrix}
  \phi_{1}\otimes \phi' &   0  &  0   &    0      \\
  0   & \phi_{1}\otimes \phi' &  0   &    0  \\
   0   &  0   & \phi_{2}\otimes \phi' &    0 \\
    0    &  0  &  0     &   \phi_{2}\otimes \phi'
\end{bmatrix}$ and $Q=\begin{bmatrix}
  \psi_{1}\otimes \psi' &   0  &  0   &    0      \\
  0   & \psi_{1}\otimes \psi' &  0   &    0  \\
   0   &  0   & \psi_{2}\otimes \psi' &    0 \\
    0    &  0  &  0     &   \psi_{2}\otimes \psi'
\end{bmatrix}$
Next, $$(X_{1}\oplus X_{2})\widetilde{\otimes} X'=((\phi_{1},\psi_{1})\oplus (\phi_{2},\psi_{2}))\widetilde{\otimes} (\phi',\psi')$$
\[=(\begin{bmatrix}
  \phi_{1} &          0      \\
    0               &   \phi_{2}
\end{bmatrix}, \begin{bmatrix}
  \psi_{1}   &          0      \\
    0               &   \psi_{2}
\end{bmatrix})\widetilde{\otimes}(\phi',\psi')\]\\
\[=(\begin{bmatrix}
  \begin{bmatrix}
  \phi_{1} &          0      \\
    0               &   \phi_{2}
\end{bmatrix}\otimes\phi' &          0      \\
    0               &   \begin{bmatrix}
  \phi_{1} &          0      \\
    0               &   \phi_{2}
\end{bmatrix}\otimes\phi'
\end{bmatrix}, \begin{bmatrix}
  \begin{bmatrix}
  \psi_{1}   &          0      \\
    0               &   \psi_{2}
\end{bmatrix}\otimes\psi'   &  0  \\
    0               &   \begin{bmatrix}
  \psi_{1}   &          0      \\
    0               &   \psi_{2}
\end{bmatrix}\otimes \psi'
\end{bmatrix})\]\\
\[=(\begin{bmatrix}
  \phi_{1}\otimes \phi' &   0  &  0   &    0      \\
  0   & \phi_{2}\otimes \phi' &  0   &    0  \\
   0   &  0   & \phi_{1}\otimes \phi' &    0 \\
    0    &  0  &  0     &   \phi_{2}\otimes \phi'
\end{bmatrix}, \begin{bmatrix}
  \psi_{1}\otimes \psi' &   0  &  0   &    0      \\
  0   & \psi_{2}\otimes \psi' &  0   &    0  \\
   0   &  0   & \psi_{1}\otimes \psi' &    0 \\
    0    &  0  &  0     &   \psi_{2}\otimes \psi'
\end{bmatrix})\,\,\,\,\,\,\,\,\cdots (\sharp')\]\\

Let $M=\begin{bmatrix}
  \phi_{1}\otimes \phi' &   0  &  0   &    0      \\
  0   & \phi_{2}\otimes \phi' &  0   &    0  \\
   0   &  0   & \phi_{1}\otimes \phi' &    0 \\
    0    &  0  &  0     &   \phi_{2}\otimes \phi'
\end{bmatrix}$ and $N=\begin{bmatrix}
  \psi_{1}\otimes \psi' &   0  &  0   &    0      \\
  0   & \psi_{2}\otimes \psi' &  0   &    0  \\
   0   &  0   & \psi_{1}\otimes \psi' &    0 \\
    0    &  0  &  0     &   \psi_{2}\otimes \psi'
\end{bmatrix}$
\\
 Recall that $X_{1}$, $X_{2}$ and $X'$ are respectively of sizes $n_{1}$, $n_{2}$ and $p$. By definition of a morphism in $MF(fg)$, we find the pair of matrices $(\delta,\beta)$ such that the following diagram commutes:
$$\xymatrix@ R=0.6in @ C=.85in{K[[x,y]]^{2(n_{1}+n_{2})p} \ar[r]^{Q} \ar[d]_{\delta} &
K[[x,y]]^{2(n_{1}+n_{2})p} \ar[d]^{\beta} \ar[r]^{P} & K[[x,y]]^{2(n_{1}+n_{2})p}\ar[d]^{\delta}\\
K[[x,y]]^{2(n_{1}+n_{2})p} \ar[r]_{N} & K[[x,y]]^{2(n_{1}+n_{2})p}\ar[r]_{M} & K[[x,y]]^{2(n_{1}+n_{2})p}}$$
We claim that if we choose $(\delta=M,\beta=P)$, then the above diagram will commute viz. $\delta P=M\beta$ and $N\delta=\beta Q$. It is clear that with our choice of $\delta$ and $\beta$, the equation
$\delta P=M\beta$ obviously holds. As for the equality $N\delta=\beta Q$, it is true just in case the following matrix equation holds:

$ \begin{bmatrix}
  \psi_{1}\otimes \psi' &   0  &  0   &    0      \\
  0   & \psi_{2}\otimes \psi' &  0   &    0  \\
   0   &  0   & \psi_{1}\otimes \psi' &    0 \\
    0    &  0  &  0     &   \psi_{2}\otimes \psi'
\end{bmatrix}\begin{bmatrix}
  \phi_{1}\otimes \phi' &   0  &  0   &    0      \\
  0   & \phi_{2}\otimes \phi' &  0   &    0  \\
   0   &  0   & \phi_{1}\otimes \phi' &    0 \\
    0    &  0  &  0     &   \phi_{2}\otimes \phi'
\end{bmatrix}$\\
$=\begin{bmatrix}
  \phi_{1}\otimes \phi' &   0  &  0   &    0      \\
  0   & \phi_{1}\otimes \phi' &  0   &    0  \\
   0   &  0   & \phi_{2}\otimes \phi' &    0 \\
    0    &  0  &  0     &   \phi_{2}\otimes \phi'
\end{bmatrix}\begin{bmatrix}
  \psi_{1}\otimes \psi' &   0  &  0   &    0      \\
  0   & \psi_{1}\otimes \psi' &  0   &    0  \\
   0   &  0   & \psi_{2}\otimes \psi' &    0 \\
    0    &  0  &  0     &   \psi_{2}\otimes \psi'
\end{bmatrix}
$\\
Hence, we need to check if each of the following smaller equalities hold.
\begin{enumerate}[a.]
  \item ($\psi_{1}\otimes \psi')(\phi_{1}\otimes \phi')=(\phi_{1}\otimes \phi')(\psi_{1}\otimes \psi')$; i.e., $\psi_{1}\phi_{1}\otimes \psi'\phi'=\phi_{1}\psi_{1}\otimes \phi'\psi'$\\
      (This holds since $\psi_{1}\phi_{1}=\phi_{1}\psi_{1}$ and $\phi'\psi'= \psi'\phi'$ because matrices involved in a matrix factorization commute by definition \ref{defn matrix facto of polyn}.
  \item ($\psi_{2}\otimes \psi')(\phi_{2}\otimes \phi')=(\phi_{1}\otimes \phi')(\psi_{1}\otimes \psi')$; i.e., $\psi_{2}\phi_{2}\otimes \psi'\phi'=\phi_{1}\psi_{1}\otimes \phi'\psi'$; i.e., $f\cdot I_{n_{2}}\otimes g\cdot I_{p}=f\cdot I_{n_{1}}\otimes g\cdot I_{p}$ since $X_{i}=(\phi_{i},\psi_{i})$, for $i=1,2$; is a matrix factorization of $f$.\\
      i.e., $fg\cdot I_{n_{2}}\otimes I_{p}=fg\cdot I_{n_{1}}\otimes I_{p}$\\
      This holds since $n_{1}=n_{2}$ as $X_{1}$ and $X_{2}$ are matrix factorizations of the same size by assumption.
  \item ($\psi_{1}\otimes \psi')(\phi_{1}\otimes \phi')=(\phi_{2}\otimes \phi')(\psi_{2}\otimes \psi')$; i.e., $\psi_{1}\phi_{1}\otimes \psi'\phi'=\phi_{2}\psi_{2}\otimes \phi'\psi'$\\
      An argument similar to the previous one shows that this identity holds.

  \item ($\psi_{2}\otimes \psi')(\phi_{2}\otimes \phi')=(\phi_{2}\otimes \phi')(\psi_{2}\otimes \psi')$; i.e., $\psi_{2}\phi_{2}\otimes \psi'\phi'=\phi_{2}\psi_{2}\otimes \phi'\psi'$\\
    A straightforward computation shows that this equality holds.
\end{enumerate}
Since the four foregoing equalities hold, the above diagram commutes.\\
The desired isomorphism now easily follows from $(\sharp)$ and $(\sharp')$. In fact,
observe that if $c_{i} \leftrightarrow c_{j}$ (resp. $r_{i}\leftrightarrow r_{j} $) stands for the operation of interchanging column $i$ and column $j$ (resp. row $i$ and row $j$) then $\sharp$ and $\sharp'$ can be obtained from one another first by applying $c_{2}\leftrightarrow c_{3}$ and then $r_{2}\leftrightarrow r_{3}$.

\item The proof showing that there is a map from $X' \widetilde{\otimes}(X_{1}\oplus X_{2})$ to $(X'\widetilde{\otimes}X_{1})\oplus (X'\widetilde{\otimes}X_{2})$ is similar to the proof we gave in part 1. The isomorphism follows from the commutativity of $\widetilde{\otimes}$.
\end{enumerate}
\end{proof}

In the next section, we give an application of the tensor products discussed in the previous sections. In fact, we will combine the tensor product of matrix factorizations $\widehat{\otimes}$ (or any of its variants) and the newly defined multiplicative tensor product of matrix factorization $\widetilde{\otimes}$ (or its variant) to improve the standard algorithm for factoring polynomials on the class of \textit{summand-reducible polynomials} defined in section \ref{chap: improvement of the std method}.

\section{The standard method improved on the class of summand-reducible polynomials} \label{chap: improvement of the std method}
In this section, we first recall a standard algorithm for factoring polynomials which dates to the 1980s when Kn\"{o}rrer exploited it to prove his celebrated periodicity theorem (cf. theorem 2.1 \cite{brown2016knorrer}). This standard technique, usually referred to as the standard method \cite{crisler2016matrix} for factoring polynomials, builds matrix factorizations of sums of polynomials from "factorizations" of their summands. One conspicuous downside of this algorithm is that for each new summand that is added to the polynomial being factorized, the size (i.e., the number of rows and columns) of the matrix factorization doubles. We define a class of polynomials and improve the standard method on this class. In fact, we define a \textit{summand-reducible polynomial} to be one that can be written in the form $f=t_{1}+\cdots + t_{s}+ g_{11}\cdots g_{1m_{1}} + \cdots + g_{l1}\cdots g_{lm_{l}}$ under some specified conditions where each $t_{k}$ is a monomial and each $g_{ji}$ is a sum of monomials. We then use tools which were not available in the days the standard method was developed, namely $\widehat{\otimes}$ (or any of its variants) and $\widetilde{\otimes}$ (or its variant), to improve the standard method for matrix factorization of polynomials on this class. We prove that if $p_{ji}$ is the number of monomials in $g_{ji}$, then there is an improved version of the standard method for factoring $f$ which produces factorizations of size $2^{\prod_{i=1}^{m_{1}}p_{1i} + \cdots + \prod_{i=1}^{m_{l}}p_{li} - (\sum_{i=1}^{m_{1}}p_{1i} + \cdots + \sum_{i=1}^{m_{l}}p_{li})}$ times smaller than the size one would normally obtain with the standard method.
It is not superfluous to mention that each time we use any of the variants of $\widehat{\otimes}$ or $\widetilde{\otimes}$, we obtain a new matrix factorization of the same size. \\
In our presentation, we limit ourselves to polynomials in the ring $S=\mathbb{R}[x_{1}\dots,x_{n}]$ where $\mathbb{R}$ is the set of real numbers. Moreover, in an attempt to simplify the presentation we will mostly only use two of the five tensor products mentioned above, namely $\widehat{\otimes}$ and $\widetilde{\otimes}$.

\subsection{The standard method for factoring polynomials} \label{sec: the std method}
\textbf{Introduction}\\
Eisenbud \cite{eisenbud1980homological} proved that using matrices, both reducible and irreducible polynomials in S can be factorized. He showed that the matrix factorizations of the polynomial $f$ are intimately related
to homological properties of modules over the quotient ring $S/(f)$, known as the hypersurface ring. \cite{knorrer1987cohen} and \cite{buchweitz1987cohen} contain more background on the connection between matrix factorizations and algebraic geometry. These papers have details on the connection that exists between matrix factorizations and maximal Cohen-Macaulay Modules. In this section, we describe a way to construct matrix factorizations of a polynomial without resorting to the homological methods that Eisenbud introduced.\\
\textbf{The standard method}\\
Here, we recall the standard technique for factoring polynomials using matrices.
\begin{proposition} \cite{crisler2016matrix}
  For $i,j\in \{1,2\}$, let $(C_{i},D_{i})$ denote an $n\times n$ matrix factorization of the polynomial $f_{i}\in S$. In addition, assume that the matrices $C_{i}$ and $D_{j}$ commute when $i\neq j.$ Then the matrices

  $$
\begin{pmatrix}
  C_{1}
  & \rvline & -D_{2} \\
\hline
  C_{2} & \rvline &
  D_{1}
\end{pmatrix},
\begin{pmatrix}
  D_{1}
  & \rvline & D_{2} \\
\hline
  -C_{2} & \rvline &
  C_{1}
\end{pmatrix}
$$
give a $2n\times 2n$ matrix factorization of $f_{1} + f_{2}$.
\end{proposition}
 The following consequence of the foregoing result is actually the basis for the main construction of the standard algorithm for matrix factorization of polynomials.
 \begin{corollary}\cite{crisler2016matrix}
 If $(C,D)$ is an $n\times n$ matrix factorization of $f$ and $g,h$ are two polynomials, then
 $$
\begin{pmatrix}
  C
  & \rvline & -gI_{n} \\
\hline
  hI_{n} & \rvline &
  D
\end{pmatrix},
\begin{pmatrix}
  D
  & \rvline & gI_{n} \\
\hline
  -hI_{n} & \rvline &
  C
\end{pmatrix}
$$

give a $2n\times 2n$ matrix factorization of $f + gh$.
 \end{corollary}
 \begin{proof}
   Since the matrices $gI_{n}$ and $hI_{n}$ commute with all $n\times n$ matrices, the proof follows from the previous proposition.
 \end{proof}

 Thanks to this corollary, one can inductively construct matrix factorizations of polynomials of the form:\\
 $$f = f_{k} = g_{1}h_{1} + g_{2}h_{2}+ \cdots + g_{k}h_{k}.$$
 For $k=1$, we have $f=g_{1}h_{1}$ and clearly $[g_{1}][h_{1}]=[g_{1}h_{1}]=[f_{1}]$ is a $1\times 1$ matrix factorization. Next, assume that $C$ and $D$ are matrix factorizations of $f_{k-1}$, i.e., $CD=If_{k-1}$ where $I$ is the identity matrix of the right size. Hence, using the foregoing corollary, we obtain a matrix factorization of $f_{k}$:
 $$
(\begin{pmatrix}
  C
  & \rvline & -g_{k}I_{n} \\
\hline
  h_{k}I_{n} & \rvline &
  D
\end{pmatrix},
\begin{pmatrix}
  D
  & \rvline & g_{k}I_{n} \\
\hline
  -h_{k}I_{n} & \rvline &
  C
\end{pmatrix})
$$
 \begin{definition} \cite{crisler2016matrix}
   The foregoing algorithm is called \textbf{the standard method} for factoring polynomials.
 \end{definition}
 Since every polynomial can be expressed as a sum of finitely many monomials, this algorithm can be used to produce matrix factorizations of any polynomial.\\
 Though this algorithm works for any polynomial, it has a conspicuous downside. The sizes of factorizations grow very quickly due to the fact that for every new summand $g_{n}h_{n}$ added to the polynomial, the factorizations double in size. It is easy to see that with this method, to factor a polynomial with $k$ summands, say $$f_{k}=g_{1}h_{1} + g_{2}h_{2}+ \cdots + g_{k}h_{k},$$ one obtains matrices of size $2^{k-1}$. This entails that these factorizations can grow extremely large very quickly. For example if $k=7$, we will obtain matrices of size $ 2^{6}= 64$ and for $k=10$, we will obtain matrices of size $2^{9} = 512$.

 Diveris and Crisler in \cite{crisler2016matrix} improved this standard algorithm on a special class of polynomials: polynomials that are sums of squares i.e., $f_{n}=x_{1}^{2} + \cdots + x_{n}^{2}$, for $n\leq 8$. The resulting factorizations they obtained have smaller matrices than one would obtain using the standard method. \cite{buchweitz1987cohen} studies matrix factorizations over quadratic hypersurfaces and also contains factorizations of $f_{n}=x_{1}^{2} + \cdots + x_{n}^{2}$. The authors of \cite{buchweitz1987cohen} first prove that there is an equivalence of categories between matrix factorizations of $f_{n}$ and graded modules over a Clifford algebra associated to $f_{n}$. They then exploit this technique to generate matrix factorizations of $f_{n}$. In \cite{crisler2016matrix} it is observed that this technique can be used to generate minimal matrix factorizations of polynomials $f_{n}$ for all $n\geq 1$. In contrast, Diveris and Crisler use an elementary approach based on matrix algebra. They remark that their algorithm produces a factorization of $f_{8}$ with just $8\times 8$ matrices whereas the standard method will produce a factorization of size $128\times 128$. Moreover, they state that the results in \cite{buchweitz1987cohen} actually prove that their factorizations of $f_{n}$ have the smallest possible size for $1\leq n \leq 8$. In fact, the authors of \cite{buchweitz1987cohen} prove that for $n\geq 8$, the smallest possible matrix factorization for $f_{n}$ is bounded below by $2^{\frac{n-2}{2}}\times 2^{\frac{n-2}{2}}$. This lower bound on the smallest matrix factorization is a crucial argument they use in deducing Hurwitz's theorem that no real composition algebra of dimension $n$ exist for $n\neq 1,2,4$ or $8$. Furthermore, they show that a necessary condition for the existence of a real composition algebra of dimension $n$ is that $f_{n}$ admits a matrix factorization of size $n\times n$. Since for $n > 8$, $n<2^{\frac{n-2}{2}}$, they deduce that no composition algebra of dimension $n$ exists when $n > 8$.  \\
 We now give some examples to illustrate the standard method.

 \begin{example} \label{exple: simple example using standard method}
   Let $h=x^{2}+ z^{2}$. Then using the standard method, a matrix factorization of $h$ is
   \[
M=
  (\begin{bmatrix}
    x & -z \\
    z & x
  \end{bmatrix},\begin{bmatrix}
    x & z \\
    -z & x
  \end{bmatrix})
\]

 \end{example}
 \begin{example} \label{exple: good matrix facto of g}
   Let $g= xy+x^{2}z+yz^{2}$. We use the standard method to find a matrix factorization of $g$.
   First a matrix factorization of $xy+x^{2}z$ is
    \[
  (\begin{bmatrix}
    x & -x^{2} \\
    z & y
  \end{bmatrix},\begin{bmatrix}
    y & x^{2} \\
    -z & x
  \end{bmatrix})
\]
Hence, a matrix factorization of $g= xy+x^{2}z+yz^{2}$ is then:
 \[
  P=(\begin{bmatrix}
    x & -x^{2} & -y & 0\\
    z &   y    &  0 & -y\\
 z^{2} &  0    &  y &  x^{2}\\
     0 &  z^{2}&  -z & x
  \end{bmatrix},\begin{bmatrix}
     y & x^{2} &  y & 0\\
    -z &   x   &  0 & y\\
 -z^{2} &  0    &  x &  -x^{2}\\
     0 &  -z^{2}&  z & y
  \end{bmatrix})
\]
\end{example}
 \begin{example} \label{exple: bad matrix facto of g}
 Let $l= xy+(x^{2}+yz)z$. Observe that $l=g$ where $g$ is given in example \ref{exple: good matrix facto of g}. We use the standard method and quickly find a matrix factorization of $l$:

  \[
  Q=(\begin{bmatrix}
    x & -(x^{2}+yz) \\
    z & y
  \end{bmatrix},\begin{bmatrix}
    y & x^{2}+yz \\
    -z & x
  \end{bmatrix})
\]
 \end{example}
We observe that the factorization we obtain for $l$ is not as nice as the one we obtain for $g$, in the sense that the complexity of some entries in the factorization of $l$ is higher than what we have for $g$.
For instance, in $Q$ the entry $-(x^{2}+yz)$ is more complex than all the entries in $P$ of example \ref{exple: good matrix facto of g}. This shows that it is better to use the expanded version of a polynomial to find its matrix factorization.\\
Here is another more involved example to further illustrate the foregoing point.
 \begin{example}
   Let $f= xy +(xy+x^{2}z+yz^{2})(x^{2}+z^{2})$. We use the standard method and quickly find a matrix factorization of $f$:

     \[
  B=(\begin{bmatrix}
    x & -(xy+x^{2}z+yz^{2}) \\
    (x^{2}+z^{2}) & y
  \end{bmatrix},\begin{bmatrix}
    y & (xy+x^{2}z+yz^{2}) \\
    -(x^{2}+z^{2}) & x
  \end{bmatrix})
\cdots  \ddag\]

 \end{example}
 Note that the matrices obtained in $B$ are not satisfactory because the complexity of $g= xy+x^{2}z+yz^{2}$ and of $h=x^{2}+z^{2}$ which are entries in $B$ could be lower had we first expanded $f$ as follows:
 $$f= xy +(xy+x^{2}z+yz^{2})(x^{2}+z^{2})=xy + xyx^{2} + x^{2}zx^{2} + yz^{2}x^{2} + xyz^{2} + x^{2}z^{3} + yz^{4}. $$
 With this expanded version of $f$, if we apply the foregoing algorithm, it is easy to see that we will obtain factorizations of $f$ which are better than the one in $\ddag$. We say better in the sense that the entries in the matrices will not be sums of monomials (e.g. $(xy+x^{2}z+yz^{2})$ in $\ddag$) but simply monomials and thus providing a more interesting factorization.
 \begin{quote}
   So, we make the important assumption that before applying the standard method to a given polynomial, it has to be written in its expanded form.
 \end{quote}
     But this comes at a price! The size of the factorizations becomes big. \\
 In fact, in order to obtain a matrix factorization of $f$ in which the entries are monomials and not sums of monomials, since $f$ has 7 summands, using the standard method would yield matrix factorizations of size $2^{7-1}=2^{6}=64$. \\
 Maybe this is due to the fact that with the standard method, even if one knows the matrix factorizations of some summands in a polynomial, there is no way to utilize them. Observe that $f= xy +(xy+x^{2}z+yz^{2})(x^{2}+z^{2})$ and from examples \ref{exple: simple example using standard method} and \ref{exple: good matrix facto of g}, we have at hand the factorizations of each factor in the product $(xy+x^{2}z+yz^{2})(x^{2}+z^{2})$. But there is no way to exploit this information using the standard algorithm.\\
  In the following section, we will exploit this information to obtain better results on the size of the matrix factorizations.

\subsection{The improved algorithm}
Here, we use the tensor product of matrix factorizations $\widehat{\otimes}$ (or any of its variants) and the multiplicative tensor product of matrix factorizations $\widetilde{\otimes}$ (or its variant) to improve the standard algorithm for factoring polynomials that are \textit{(simple) summand-reducible} (cf. definition \ref{defn simple summand-reducible polynomial} \& \ref{defn summand reducible polynomials}). In fact, we show that our approach produces factorizations that are of smaller sizes than the factorizations produced by the standard method.
We first define a class of polynomials which are such that they can be rewritten with less summands than what appears in their expanded form.
\begin{definition} \label{defn simple summand-reducible polynomial}
  A polynomial $f$ is said to be \textbf{simple summand-reducible} if it can be written in the form
  $$f= g_{1}h_{1} + g_{2}h_{2}+ \cdots + g_{k}h_{k},\,\,k\geq 2$$ where:
  \begin{enumerate}
  \item For at least one $i\in \{1,2, \cdots, k\}$, $g_{i}h_{i}$ is the product of two polynomials ($g_{i}$ and $h_{i}$), s.t. if $n_{g_{i}}=$ number of monomials in $g_{i}$ and $n_{h_{i}}=$ number of monomials in $h_{i}$, then $n_{g_{i}}n_{h_{i}}\geq 6$.
  \item If there is an $i\in \{1,2, \cdots, k\}$ s.t. the expanded form of $g_{i}h_{i}$ belongs to $ \{a^{n}-b^{n}, a^{n}+b^{n}, (a^{2} \pm 2ab + b^{2});  n\in \mathbb{N}-\{0\}\}$, (where $a$ and $b$ could be single variables or product of variables), then instead expand $g_{i}h_{i}$ to have a smaller expression.
  \end{enumerate}
\end{definition}
The reason for denying that the two factors in $g_{i}h_{i}$ should not be simultaneously $(a \pm b)$ is that
when expanded, $(a \pm b)^{2}$ has three summands; thus, the standard method would yield a factorization of size $2^{3-1}=4$ and as we shall see (cf. theorem \ref{thm improved standard method}), our new algorithm will produce factorizations of size $2(2^{1})(2^{1})=8$. Indeed, we see that there will be no profit in having the products $(a \pm b)^{2}$ in our foregoing definition.
Therefore, it is better to write the expression that has three summands. The explanation of the other part of the second point in the definition is that factorizing $ \{a^{n}-b^{n}, a^{n}+b^{n}\}$ will instead cause the size of the matrices to grow big as we will have more monomials after factorizing.\\
The idea is that we want to be able to write polynomials with less summands than the ones that appear in their expanded form. Because, the fewer the number of summands, the smaller the size of the matrix factorizations as we will soon see with the help of the tensor products $\widehat{\otimes}$ (or any of its variants) and $\widetilde{\otimes}$ (or its variant). \\
In the remainder of this section, we will only use the operations $\widehat{\otimes}$ and $\widetilde{\otimes}$ without explicitly alluding to their variants though it is clear that in place of each of these operations, any of their variants could be used to produce the desired results.
\\ Consider the following polynomials for which the first condition of definition \ref{defn simple summand-reducible polynomial} fails:
\begin{enumerate}
  \item Let $n,m\in \mathbb{N}$, $x^{n}-y^{m}$ is not simple summand-reducible.
  \item $xy+(x-y)(x^{2} + 2xy + y^{2})= xy+x^{3}-y^{3}$ is not simple summand-reducible. Note that it is better to write $g_{i}h_{i}=(x-y)(x^{2} + 2xy + y^{2})$ as $x^{3}-y^{3}$ because the latter expression produces smaller factorizations than the former (as will be made clear in the proof of theorem \ref{thm improved standard method}).
      This example explains why we need the second condition in definition \ref{defn simple summand-reducible polynomial}.
 \item $x^{2}y + x^{2}z^{3}$ is not simple summand-reducible. Note that it is the expanded form of $x^{2}(y + z^{3})$ which is a $1\times 1$ matrix factorization of $x^{2}y + x^{2}z^{3}$.
\end{enumerate}
In contrast:
\begin{example} \label{exples of simple summand-red polyn}
  The following are simple summand-reducible polynomials:
  \begin{enumerate}
    \item $x^{2}+ yx^{3}+ zx^{4} + yz^{2}x^{2} + xy^{3} + x^{2}zy^{2} + z^{2}y^{3}= x^{2}+ (xy + x^{2}z + yz^{2})(x^{2} + y^{2})$.
    \item  $x^{3}-y^{3}+ (xy + yz^{2})(x^{2} + y^{2} + z)$
    \item $f=x^{2} + xyz + yx^{3}+ zx^{4} + yz^{2}x^{2} + xy^{3} + x^{2}zy^{2} + z^{2}y^{3} + xyz + x^{2}z^{2} + yz^{3}=x^{2} + xyz + (xy + x^{2}z + yz^{2})(x^{2} + y^{2} + z)$.
  \end{enumerate}
\end{example}
\begin{remark}
  Observe that our definition \ref{defn simple summand-reducible polynomial} mostly targets polynomials with more than six monomials because factorizations obtained with the standard method begin to be of considerable sizes such that it is even impossible to properly put some of them on an $A_{4}$ sheet of paper.
\end{remark}
\begin{definition}
  A polynomial $f$ is said to be \textbf{simple summand-reduced} if it is in the form $f = g_{1}h_{1} + g_{2}h_{2} + \cdots + g_{k}h_{k}$ described in definition \ref{defn simple summand-reducible polynomial}.
\end{definition}

We now want to derive a procedure for finding matrix factors of simple summand-reduced polynomials that are of smaller sizes than the size obtained with the standard method.\\
 
\textbf{Nota Bene:} With the standard method, even if one knows matrix factorizations of polynomials $f$ and $g$, one cannot derive from it a matrix factorization of the product $fg$ or the sum $f+g$.\\
Recall that the multiplicative tensor product of matrix factorization $\widetilde{\otimes}$ developed in this paper produces a matrix factorization of the product of two polynomials from the matrix factorizations of each of these polynomials.
In theorems \ref{thm improved standard method} and \ref{thm improved algo for summand-red polyn}, $\widetilde{\otimes}$ is one of the crucial ingredients in the improved algorithm used to obtain better matrix factorizations of polynomials in the sense that, the size of the matrix factorizations we obtain is smaller than the size one would normally obtain with the standard method. In fact, $\widetilde{\otimes}$ helps in reducing the size of the factorizations as we shall notice in the proofs of theorems \ref{thm improved standard method} and \ref{thm improved algo for summand-red polyn}. We will have a foretaste of the reduction power of $\widetilde{\otimes}$ in action in Example \ref{exple: example before the first thm} with the polynomial $(xy+x^{2}z+yz^{2})(x^{2}+z^{2})$. In fact, using $\widetilde{\otimes}$, we obtain a matrix factorization of size $2(4)(2)=16$ but if we expand it and use the standard method, we would have 6 monomials and consequently a matrix factorization of size $2^{6-1}=32$, this is twice the size we obtain with $\widetilde{\otimes}$.\\
Another crucial ingredient used in the improved algorithm is Yoshino's tensor product of matrix factorizations $\widehat{\otimes}$, because it produces a matrix factorization of the sum of two polynomials from the matrix factorizations of each of these polynomials.\\

The following example gives an idea on how to use the operations $\widehat{\otimes}$ and $\widetilde{\otimes}$ to find the sizes of matrix factors of simple summand-reduced polynomials.

\begin{example} \label{exple: example before the first thm}
  Compare the size of a matrix factorization of $f= xy +(xy+x^{2}z+yz^{2})(x^{2}+z^{2})$ obtained using the standard method with the size of a matrix factorization of $f$ obtained using the operations $\widehat{\otimes}$ and $\widetilde{\otimes}$.\\
  We know that $f$ in expanded form has $7$ monomials and consequently the size of a matrix factorization of $f$ obtained using the standard method will be $2^{7-1}=2^{6}=64$.
  \\
  We also know that we can use $\widetilde{\otimes}$ to find a matrix factorization of the product $(xy+x^{2}z+yz^{2})(x^{2}+z^{2})$ and then use $\widehat{\otimes}$ to obtain a matrix factorization of the sum $f= xy +(xy+x^{2}z+yz^{2})(x^{2}+z^{2})$. 
  We already found $M$ and $P$ above which are respectively matrix factorizations of $(x^{2}+z^{2})$ and $(xy+x^{2}z+yz^{2})$. $P \widetilde{\otimes} M$ will be of size $2pm=2(4)(2)=16$ where $p=4$ is the size of $P$ and $m=2$ is the size of $M$. Now, let $Q=([x],[y])$ be a $1\times 1$ matrix factorization of $xy$.
  Then, if $L= P \widetilde{\otimes} M$; $Q\widehat{\otimes} L$ would be a matrix factorization of $f= xy +(xy+x^{2}z+yz^{2})(x^{2}+z^{2})$ of size $2ql=2(1)(16)=32$, where $q=1$ is the size of the matrix factorization $Q$ and $l=16$ is the size of the matrix factorization $L$. \\
  So, the size of the factorization we obtain for $f$ using the operations $\widehat{\otimes}$ and $\widetilde{\otimes}$ is $32=\frac{64}{2}$. That is one-half the size we obtain using the standard method!  \\ 
  Finding a matrix factorization of $f$ is now easy after all the foregoing explanation. This will be done below in example \ref{exple: matrix factors presented, part I}.
\end{example}
We can now state and prove the following theorem:
\begin{theorem} \label{thm improved standard method}
There is an improved version of the standard method for factoring simple summand-reducible polynomials which produces factorizations which are at most one-half the size of the matrix factorizations one would obtain with the standard method.
\end{theorem}

\begin{proof}
  First, we construct the algorithm, then we prove that the size of the resulting matrix factorizations (for simple summand-reducible polynomials) is at most one-half the size one would obtain with the standard method.\\
  We inductively construct the matrix factorizations of simple summand-reduced polynomials using tools ($\widehat{\otimes}$ and $\widetilde{\otimes}$) that were not existing in the 1980s when the standard method was constructed.\\
  Let $f=f_{k}= g_{1}h_{1} + g_{2}h_{2}+ \cdots + g_{k}h_{k}$ be a simple summand-reduced polynomial. Without loss of generality, we consider the last summand to be a product of sums of polynomials.\\
  For $k=2$: \\

  Since $g_{2}h_{2}$ is a product of sums of monomials, then:\\
  \begin{enumerate}
    \item Use the standard method to find a matrix factorization of $g_{2}$, call it $X_{g_{2}}=(A_{g_{2}},B_{g_{2}})$.
    \item Use the standard method to find a matrix factorization of $h_{2}$, call it $X_{h_{2}}=(A_{h_{2}},B_{h_{2}})$.
    \item Use the multiplicative tensor product $\widetilde{\otimes}$ to find $X_{g_{2}} \widetilde{\otimes} X_{h_{2}}$ which is a matrix factorization of $g_{2}h_{2}$.
  \end{enumerate}
  Next, \\
 $\bullet$ if $g_{1}h_{1}$ is a product of sums of monomials, then proceed as above to obtain a matrix factorization of $g_{1}h_{1}$ and then use $\widehat{\otimes}$ to obtain a matrix factorization of $f$.\\
$\bullet$ Else $g_{1}h_{1}$ is not a product of sums of monomials, thus $[g_{1}][h_{1}]$ is a $1\times 1$ matrix factorization of $g_{1}h_{1}$. Hence, use $\widehat{\otimes}$ to find a matrix factorization of $f$.\\

 Next, assume $(A,B)$ is a matrix factorization of $f_{k-1}$:\\
 \begin{enumerate}[a.]
   \item Use the standard method to find a matrix factorization of $g_{k}$,
   \item Use the standard method to find a matrix factorization of $h_{k}$,
   \item Use the multiplicative tensor product $\widetilde{\otimes}$ to find a matrix factorization of $g_{k}h_{k}$, call it $(C,D)$.
   \item Finally, use the tensor product of matrix factorizations $\widehat{\otimes}$ to find $(A,B) \widehat{\otimes} (C,D)$ which is a matrix factorization of $f=f_{k-1} + g_{k}h_{k}$.
 \end{enumerate}
We now show that the size of the factorizations we obtain using the improved method is at most one-half the size one would obtain using the standard method.\\
Let $f$ be a simple summand-reducible polynomial with $n$ monomials in its expanded form. Then the standard method tells us that matrix factorizations of $f$ will be of size $2^{n-1}$.\\
Let $f=f_{k}= g_{1}h_{1} + g_{2}h_{2}+ \cdots + g_{k}h_{k}$ be the simple summand-reduced form of $f$. Then at least one of the $g_{i}h_{i}$ should be the product of sums of monomials. Assume $f$ has just one such summand. Without loss of generality (WLOG), let it be $g_{k}h_{k}$. We know by definition of simple summand-reducible polynomial that the number of monomials in $g_{k}$ times the number of monomials in $h_{k}$ should be at least $6$. Suppose it is $6$ (the reason why this assumption is sufficient is found in theorem \ref{Cor. improved algo on t+gh} and the comment that follows it.) and suppose WLOG that $g_{k}$ and $h_{k}$ are respectively the sums of three  and two monomials. Next, using the standard method, we can find factorizations of $g_{k}$ and $h_{k}$ respectively of sizes $2^{3-1}=4$ and $2^{2-1}=2$. Hence, using the multiplicative tensor product $(\widetilde{\otimes})$ of matrix factorizations, we find a matrix factorization of $g_{k}h_{k}$ of size $2(4)(2)=16=2^{4}$. \\
Now $f$ had $n$ monomials in its expanded form, meaning that if we take away the $6$ monomials that made up $g_{k}h_{k}$, we will be left with $n-6$ monomials in $f$. Thus the standard method would yield a matrix factorization of size $2^{(n-6)-1}$ of $f-g_{k}h_{k}$. Finally, we use the tensor product of matrix factorizations $(\widehat{\otimes})$ to find a matrix factorization of $f=(f-g_{k}h_{k}) + g_{k}h_{k}$ since we have factorizations of each of these two summands. Thus, we obtain a factorization of size $2(2^{(n-6)-1})(2^{4})=2^{n-2}$ and this is one-half the size we obtained when we solely used the standard method.
\end{proof}

\begin{example} \label{exple: matrix factors presented, part I}
Use the improved algorithm to factorize the polynomial \\$f= xy +(x^{2}+z^{2})(xy+x^{2}z+yz^{2})$ given in Example \ref{exple: example before the first thm}. \\
$f=xy + hg$ where $h=x^{2}+z^{2}$ and $g=xy+x^{2}z+yz^{2}$. In examples \ref{exple: simple example using standard method} and \ref{exple: good matrix facto of g}, we used the standard method to find matrix factorizations of the polynomials $h$ and $g$ which are respectively:
\[
M=(\phi_{h},\psi_{h})=
  (\begin{bmatrix}
    x & -z \\
    z & x
  \end{bmatrix},\begin{bmatrix}
    x & z \\
    -z & x
  \end{bmatrix})
\]
and 
\[
  P=(\phi_{g},\psi_{g})=(\begin{bmatrix}
    x & -x^{2} & -y & 0\\
    z &   y    &  0 & -y\\
 z^{2} &  0    &  y &  x^{2}\\
     0 &  z^{2}&  -z & x
  \end{bmatrix},\begin{bmatrix}
     y & x^{2} &  y & 0\\
    -z &   x   &  0 & y\\
 -z^{2} &  0    &  x &  -x^{2}\\
     0 &  -z^{2}&  z & y
  \end{bmatrix})
\]
According to the proof of Theorem \ref{thm improved standard method}, to find a matrix factorization for $f$, we need to:
\begin{enumerate}
  \item First of all find a matrix factorization of the product $hg$ using the multiplicative  tensor product $\widetilde{\otimes}$ as follows:\\
$M \widetilde{\otimes} P = (\phi_{h},\psi_{h})\widetilde{\otimes}(\phi_{g},\psi_{g})=(\phi_{hg},\psi_{hg})$ where $$(\phi_{hg},\psi_{hg})= (\begin{bmatrix}
    \phi_{h}\otimes \phi_{g} & 0 \\
    0             &   \phi_{h}\otimes \phi_{g}
  \end{bmatrix},\begin{bmatrix}
   \psi_{h}\otimes \psi_{g} & 0 \\
    0             &   \psi_{h}\otimes \psi_{g}
  \end{bmatrix})$$
with \\
$\,\,\,\,\,\,\,\,\,\,\phi_{h}\otimes \phi_{g}\,\,\,\,\,\,\,\,\,\,\,\,\,\,\,\,\,\,= 
\begin{bmatrix}
    x & -z \\
    z & x
  \end{bmatrix} \otimes \begin{bmatrix}
    x & -x^{2} & -y & 0\\
    z &   y    &  0 & -y\\
 z^{2} &  0    &  y &  x^{2}\\
     0 &  z^{2}&  -z & x
  \end{bmatrix}$
$$=\begin{bmatrix}
    x\begin{bmatrix}
    x & -x^{2} & -y & 0\\
    z &   y    &  0 & -y\\
 z^{2} &  0    &  y &  x^{2}\\
     0 &  z^{2}&  -z & x
  \end{bmatrix} & -z\begin{bmatrix}
    x & -x^{2} & -y & 0\\
    z &   y    &  0 & -y\\
 z^{2} &  0    &  y &  x^{2}\\
     0 &  z^{2}&  -z & x
  \end{bmatrix} \\
    z\begin{bmatrix}
    x & -x^{2} & -y & 0\\
    z &   y    &  0 & -y\\
 z^{2} &  0    &  y &  x^{2}\\
     0 &  z^{2}&  -z & x
  \end{bmatrix} & x\begin{bmatrix}
    x & -x^{2} & -y & 0\\
    z &   y    &  0 & -y\\
 z^{2} &  0    &  y &  x^{2}\\
     0 &  z^{2}&  -z & x
  \end{bmatrix}
  \end{bmatrix} $$
  $$=\begin{bmatrix}
    x^{2} & -x^{3} & -xy & 0 & -zx & zx^{2} & zy & 0\\
    xz &   xy    &  0 & -xy & -z^{2} & -zy & 0 & zy \\
    xz^{2} &  0 &  xy &  x^{3} & -z^{3} & 0 & -zy & -zx^{2}\\
     0 &  xz^{2}&  -xz & x^{2} & 0 & -z^{3} & z^{2} & -zx\\
     zx & -zx^{2} & -zy & 0 & x^{2} & -x^{3} & -xy & 0\\
     z^{2} & zy & 0 & -zy & xz & xy & 0 & -xy \\
     z^{3} & 0 & zy & zx^{2} & xz^{2} & 0 & xy & x^{3} \\
      0 & z^{3} & -z^{2} & zx &  0 &  xz^{2}&  -xz & x^{2}
  \end{bmatrix}$$
  Thus,
   
 \begin{gather*}
  \setlength{\arraycolsep}{1.0\arraycolsep}
  \renewcommand{\arraystretch}{1.5}
  \text{\footnotesize$\displaystyle
    \phi_{hg}=\begin{pmatrix}
      x^{2} & -x^{3} & -xy & 0 & -zx & zx^{2} & zy & 0 & 0& 0& 0 & 0 & 0& 0& 0&0 \\
    xz &   xy    &  0 & -xy & -z^{2} & -zy & 0 & zy & 0& 0& 0 & 0 & 0& 0& 0&0 \\
    xz^{2} &  0 &  xy &  x^{3} & -z^{3} & 0 & -zy & -zx^{2} & 0& 0& 0 & 0 & 0& 0& 0&0 \\
     0 &  xz^{2}&  -xz & x^{2} & 0 & -z^{3} & z^{2} & -zx & 0& 0& 0 & 0 & 0& 0& 0&0 \\
     zx & -zx^{2} & -zy & 0 & x^{2} & -x^{3} & -xy & 0 & 0& 0& 0 & 0 & 0& 0& 0&0 \\
     z^{2} & zy & 0 & -zy & xz & xy & 0 & -xy & 0& 0& 0 & 0 & 0& 0& 0&0 \\
     z^{3} & 0 & zy & zx^{2} & xz^{2} & 0 & xy & x^{3} & 0& 0& 0 & 0 & 0& 0& 0&0 \\
      0 & z^{3} & -z^{2} & zx &  0 &  xz^{2}&  -xz & x^{2} & 0& 0& 0 & 0 & 0& 0& 0&0 \\
      0& 0& 0 & 0 & 0& 0& 0&0 & x^{2} & -x^{3} & -xy & 0 & -zx & zx^{2} & zy & 0\\
     0& 0& 0 & 0 & 0& 0& 0&0 & xz &   xy    &  0 & -xy & -z^{2} & -zy & 0 & zy \\
     0& 0& 0 & 0 & 0& 0& 0&0 & xz^{2} &  0 &  xy &  x^{3} & -z^{3} & 0 & -zy & -zx^{2}\\
     0& 0& 0 & 0 & 0& 0& 0&0 & 0 &  xz^{2}&  -xz & x^{2} & 0 & -z^{3} & z^{2} & -zx\\
     0& 0& 0 & 0 & 0& 0& 0&0 & zx & -zx^{2} & -zy & 0 & x^{2} & -x^{3} & -xy & 0\\
     0& 0& 0 & 0 & 0& 0& 0&0 &z^{2} & zy & 0 & -zy & xz & xy & 0 & -xy \\
     0& 0& 0 & 0 & 0& 0& 0&0 & z^{3} & 0 & zy & zx^{2} & xz^{2} & 0 & xy & x^{3} \\
      0& 0& 0 & 0 & 0& 0& 0&0 &0 & z^{3} & -z^{2} & zx &  0 &  xz^{2}&  -xz & x^{2}
 \end{pmatrix}
  $}
\end{gather*}

We now compute the second matrix factor of the product $gh$. We know that $M \widetilde{\otimes} P = (\phi_{h},\psi_{h})\widetilde{\otimes}(\phi_{g},\psi_{g})=(\phi_{hg},\psi_{hg})$. So, what we want to calculate is $\psi_{hg}$. Note that by Lemma \ref{lemma size of X tensor Y}, $M \widetilde{\otimes} P$ is of size $2(2)(4)=16$, since $(\phi_{h},\psi_{h})$ is of size $2$ and $(\phi_{g},\psi_{g})$ is of size $4$. First, we calculate $\psi_{h}\otimes \psi_{g}$:\\
$\,\,\,\,\,\,\,\,\,\,\psi_{h}\otimes \psi_{g}\,\,\,\,\,\,\,\,\,\,\,\,\,\,\,\,\,\,= 
\begin{bmatrix}
    x & z \\
    -z & x
  \end{bmatrix} \otimes \begin{bmatrix}
    y & x^{2} & y & 0\\
    -z &   x    &  0 & y\\
   -z^{2} &  0    &  x &  -x^{2}\\
     0 &  -z^{2}&  z & y
  \end{bmatrix}$
 $$=\begin{bmatrix}
    x\begin{bmatrix}
    y & x^{2} & y & 0\\
    -z &   x    &  0 & y\\
   -z^{2} &  0    &  x &  -x^{2}\\
     0 &  -z^{2}&  z & y
  \end{bmatrix} & z\begin{bmatrix}
    y & x^{2} & y & 0\\
    -z &   x    &  0 & y\\
   -z^{2} &  0    &  x &  -x^{2}\\
     0 &  -z^{2}&  z & y
  \end{bmatrix} \\
    -z\begin{bmatrix}
    y & x^{2} & y & 0\\
    -z &   x    &  0 & y\\
   -z^{2} &  0    &  x &  -x^{2}\\
     0 &  -z^{2}&  z & y
  \end{bmatrix} & x\begin{bmatrix}
    y & x^{2} & y & 0\\
    -z &   x    &  0 & y\\
   -z^{2} &  0    &  x &  -x^{2}\\
     0 &  -z^{2}&  z & y
  \end{bmatrix}
  \end{bmatrix}$$
  $$=\begin{bmatrix}
    xy & x^{3} & xy & 0 & zy & zx^{2} & zy & 0\\
    -xz &   x^{2}    &  0 & xy & -z^{2} & zx & 0 & zy \\
   -xz^{2} &  0    &  x^{2} &  -x^{3} & -z^{3} & 0 & zx & -zx^{2}\\
     0 &  -xz^{2}&  xz & xy & 0 & -z^{3} & z^{2} & zy\\
    -zy & -zx^{2} & -zy & 0 & xy & x^{3} & xy & 0\\
    z^{2} & -zx & 0 & -zy & -xz &   x^{2}    &  0 & xy\\
    z^{3} & 0 & -zx & zx^{2} & -xz^{2} &  0  &  x^{2} &  -x^{3}\\
     0 & z^{3} & -z^{2} & -zy& 0 &  -xz^{2}&  xz & xy
  \end{bmatrix}$$

Thus,

 \begin{gather*}
  \setlength{\arraycolsep}{1.0\arraycolsep}
  \renewcommand{\arraystretch}{1.5}
  \text{\footnotesize$\displaystyle
    \psi_{hg}=\begin{pmatrix}
      xy & x^{3} & xy & 0 & zy & zx^{2} & zy & 0 & 0& 0& 0 & 0 & 0& 0& 0&0\\
    -xz &   x^{2}    &  0 & xy & -z^{2} & zx & 0 & zy & 0& 0& 0 & 0 & 0& 0& 0&0\\
   -xz^{2} &  0    &  x^{2} &  -x^{3} & -z^{3} & 0 & zx & -zx^{2} & 0& 0& 0 & 0 & 0& 0& 0&0\\
     0 &  -xz^{2}&  xz & xy & 0 & -z^{3} & z^{2} & zy & 0& 0& 0 & 0 & 0& 0& 0&0\\
    -zy & -zx^{2} & -zy & 0 & xy & x^{3} & xy & 0 & 0& 0& 0 & 0 & 0& 0& 0&0\\
    z^{2} & -zx & 0 & -zy & -xz &   x^{2}    &  0 & xy & 0& 0& 0 & 0 & 0& 0& 0&0\\
    z^{3} & 0 & -zx & zx^{2} & -xz^{2} &  0  &  x^{2} &  -x^{3} & 0& 0& 0 & 0 & 0& 0& 0&0\\
     0 & z^{3} & -z^{2} & -zy& 0 &  -xz^{2}&  xz & xy & 0& 0& 0 & 0 & 0& 0& 0&0\\
      0& 0& 0 & 0 & 0& 0& 0&0 & xy & x^{3} & xy & 0 & zy & zx^{2} & zy & 0\\
      0& 0& 0 & 0 & 0& 0& 0&0 & -xz &   x^{2}    &  0 & xy & -z^{2} & zx & 0 & zy \\
      0& 0& 0 & 0 & 0& 0& 0&0 &-xz^{2} &  0    &  x^{2} &  -x^{3} & -z^{3} & 0 & zx & -zx^{2}\\
       0& 0& 0 & 0 & 0& 0& 0&0 &0 &  -xz^{2}&  xz & xy & 0 & -z^{3} & z^{2} & zy\\
      0& 0& 0 & 0 & 0& 0& 0&0 &-zy & -zx^{2} & -zy & 0 & xy & x^{3} & xy & 0\\
     0& 0& 0 & 0 & 0& 0& 0&0 & z^{2} & -zx & 0 & -zy & -xz &   x^{2}    &  0 & xy\\
      0& 0& 0 & 0 & 0& 0& 0&0 &z^{3} & 0 & -zx & zx^{2} & -xz^{2} &  0  &  x^{2} &  -x^{3}\\
      0& 0& 0 & 0 & 0& 0& 0&0 & 0 & z^{3} & -z^{2} & -zy& 0 &  -xz^{2}&  xz & xy
    \end{pmatrix}
  $}
\end{gather*}

 \item Next, from the algorithm given in Theorem \ref{thm improved standard method}, we now need to find a matrix factorization of $r=xy$ (which is the first summand in $f$). Obviously, $Q=(\phi_{r},\psi_{r})=([x],[y])$ is a matrix factorization of $xy$.
\item Finally, from our algorithm, we find a matrix factorization of $f$ by computing $Q\widehat{\otimes} (M\widetilde{\otimes} P)$ which will be of size $2(1)(16)=32$ by Lemma 2.1 of \cite{fomatati2019multiplicative} since $Q$ is of size $1$ and $(M\widetilde{\otimes} P)$ is of size $16$.\\
    In the sequel, recall that $1_{n}$ is the identity $n\times n$ matrix. We have:
    \begin{align*}
    Q\widehat{\otimes} (M\widetilde{\otimes} P) &= (\phi_{r},\psi_{r})\widehat{\otimes} (\phi_{hg}, \psi_{hg})\\
&=( \begin{bmatrix}
    \phi_{r}\otimes 1_{16}  &  1_{1}\otimes \phi_{hg}      \\
   -1_{1}\otimes \psi_{hg}  &  \psi_{r}\otimes 1_{16}
\end{bmatrix},
\begin{bmatrix}
    \psi_{r}\otimes 1_{16}  &  -1_{1}\otimes \phi_{hg}     \\
    1_{1}\otimes \psi_{hg}  &  \phi_{r}\otimes 1_{16}
\end{bmatrix}
)\\
&=( \begin{bmatrix}
    x\otimes 1_{16}  &  1\otimes \phi_{hg}      \\
   -1\otimes \psi_{hg}  &  y\otimes 1_{16}
\end{bmatrix},
\begin{bmatrix}
    y\otimes 1_{16}  &  -1\otimes \phi_{hg}     \\
    1\otimes \psi_{hg}  &  x\otimes 1_{16}
\end{bmatrix}
)\\
&=( \begin{bmatrix}
    x\otimes 1_{16}  &   \phi_{hg}      \\
   - \psi_{hg}  &  y\otimes 1_{16}
\end{bmatrix},
\begin{bmatrix}
    y\otimes 1_{16}  &  -\phi_{hg}     \\
    \psi_{hg}  &  x\otimes 1_{16}
\end{bmatrix}
)\\
&=(\phi_{rhg},\psi_{rhg})
\end{align*}

Where:\\
$\bullet$ $x\otimes 1_{16}$ (respectively $y\otimes 1_{16}$) is a $16\times 16$ diagonal matrix with $x$ (respectively $y$) on its entire diagonal. \\
$\bullet$ $\phi_{hg}$ and $\psi_{hg}$ were computed above. \\
Thus, we found a $32\times 32$ matrix factorization of $f$ viz. a matrix factorization of $f$ of size 32. Note that the standard method yields matrix factors of size $2^{7-1}=2^{6}=64=32\times 2 $. This is twice the size we obtain with the improved algorithm. \\
The matrix factors displayed below are respectively $\phi_{rhg}$ and $\psi_{rhg}$. Each of them is of size 32.

\begin{gather*}
  \setlength{\arraycolsep}{.60\arraycolsep}
  \renewcommand{\arraystretch}{1.5}
  \text{\footnotesize$\displaystyle
    \begin{pmatrix}
     x & 0& 0& 0 & 0 & 0& 0& 0& 0 & 0& 0& 0 & 0 & 0& 0& 0& x^{2} & -x^{3} & -xy & 0 & -zx & zx^{2} & zy & 0 & 0& 0& 0 & 0 & 0& 0& 0&0 \\
     0& x & 0& 0 & 0 & 0& 0& 0& 0 & 0& 0& 0 & 0 & 0& 0& 0&xz &   xy    &  0 & -xy & -z^{2} & -zy & 0 & zy & 0& 0& 0 & 0 & 0& 0& 0&0 \\
     0& 0& x & 0 & 0 & 0& 0& 0& 0 & 0& 0& 0 & 0 & 0& 0& 0& xz^{2} &  0 &  xy &  x^{3} & -z^{3} & 0 & -zy & -zx^{2} & 0& 0& 0 & 0 & 0& 0& 0&0 \\
      0& 0& 0& x & 0 & 0& 0& 0& 0 & 0& 0& 0 & 0 & 0& 0& 0& 0 &  xz^{2}&  -xz & x^{2} & 0 & -z^{3} & z^{2} & -zx & 0& 0& 0 & 0 & 0& 0& 0&0 \\
      0& 0& 0& 0 & x & 0& 0& 0& 0 & 0& 0& 0 & 0 & 0& 0& 0& zx & -zx^{2} & -zy & 0 & x^{2} & -x^{3} & -xy & 0 & 0& 0& 0 & 0 & 0& 0& 0&0 \\
      0& 0& 0& 0 & 0 & x & 0& 0& 0 & 0& 0& 0 & 0 & 0& 0& 0&z^{2} & zy & 0 & -zy & xz & xy & 0 & -xy & 0& 0& 0 & 0 & 0& 0& 0&0 \\
      0& 0& 0& 0 & 0 & 0& x& 0& 0 & 0& 0& 0 & 0 & 0& 0& 0&z^{3} & 0 & zy & zx^{2} & xz^{2} & 0 & xy & x^{3} & 0& 0& 0 & 0 & 0& 0& 0&0 \\
       0& 0& 0& 0 & 0 & 0& 0& x& 0 & 0& 0& 0 & 0 & 0& 0& 0& 0 & z^{3} & -z^{2} & zx &  0 &  xz^{2}&  -xz & x^{2} & 0& 0& 0 & 0 & 0& 0& 0&0 \\
       0& 0& 0& 0 & 0 & 0& 0& 0& x & 0& 0& 0 & 0 & 0& 0& 0& 0& 0& 0 & 0 & 0& 0& 0&0 & x^{2} & -x^{3} & -xy & 0 & -zx & zx^{2} & zy & 0\\
      0& 0& 0& 0 & 0 & 0& 0& 0& 0 & x & 0& 0 & 0 & 0& 0& 0& 0& 0& 0 & 0 & 0& 0& 0&0 & xz &   xy    &  0 & -xy & -z^{2} & -zy & 0 & zy \\
      0& 0& 0& 0 & 0 & 0& 0& 0& 0 & 0& x & 0 & 0 & 0& 0& 0&0& 0& 0 & 0 & 0& 0& 0&0 & xz^{2} &  0 &  xy &  x^{3} & -z^{3} & 0 & -zy & -zx^{2}\\
      0& 0& 0& 0 & 0 & 0& 0& 0& 0 & 0& 0& x & 0 & 0& 0& 0&0& 0& 0 & 0 & 0& 0& 0&0 & 0 &  xz^{2}&  -xz & x^{2} & 0 & -z^{3} & z^{2} & -zx\\
      0& 0& 0& 0 & 0 & 0& 0& 0& 0 & 0& 0& 0 & x & 0& 0& 0& 0& 0& 0 & 0 & 0& 0& 0&0 & zx & -zx^{2} & -zy & 0 & x^{2} & -x^{3} & -xy & 0\\
      0& 0& 0& 0 & 0 & 0& 0& 0& 0 & 0& 0& 0 & 0 & x & 0& 0& 0& 0& 0 & 0 & 0& 0& 0&0 &z^{2} & zy & 0 & -zy & xz & xy & 0 & -xy \\
      0& 0& 0& 0 & 0 & 0& 0& 0& 0 & 0& 0& 0 & 0 & 0& x & 0& 0& 0& 0 & 0 & 0& 0& 0&0 & z^{3} & 0 & zy & zx^{2} & xz^{2} & 0 & xy & x^{3} \\
       0& 0& 0& 0 & 0 & 0& 0& 0& 0 & 0& 0& 0 & 0 & 0& 0& x& 0& 0& 0 & 0 & 0& 0& 0&0 &0 & z^{3} & -z^{2} & zx &  0 &  xz^{2}&  -xz & x^{2}\\
     -xy & -x^{3} & -xy & 0 & -zy & -zx^{2} & -zy & 0 & 0& 0& 0 & 0 & 0& 0& 0&0& y & 0& 0& 0 & 0 & 0& 0& 0& 0 & 0& 0& 0 & 0 & 0& 0& 0\\
    xz &  -x^{2}    &  0 & -xy & z^{2} & -zx & 0 & -zy & 0& 0& 0 & 0 & 0& 0& 0&0&0& y& 0& 0 & 0 & 0& 0& 0& 0 & 0& 0& 0 & 0 & 0& 0& 0\\
    xz^{2} &  0    &  -x^{2} &  x^{3} & z^{3} & 0 & -zx & zx^{2} & 0& 0& 0 & 0 & 0& 0& 0&0&0& 0& y& 0 & 0 & 0& 0& 0& 0 & 0& 0& 0 & 0 & 0& 0& 0\\
     0 &  xz^{2}&  -xz & -xy & 0 & z^{3} & -z^{2} & -zy & 0& 0& 0 & 0 & 0& 0& 0&0&0& 0& 0& y & 0 & 0& 0& 0& 0 & 0& 0& 0 & 0 & 0& 0& 0\\
    zy & zx^{2} & zy & 0 & -xy & -x^{3} & -xy & 0 & 0& 0& 0 & 0 & 0& 0& 0&0&0& 0& 0& 0 & y & 0& 0& 0& 0 & 0& 0& 0 & 0 & 0& 0& 0\\
    -z^{2} & zx & 0 & zy & xz &   -x^{2}    &  0 & -xy & 0& 0& 0 & 0 & 0& 0& 0&0&0& 0& 0& 0 & 0 & y& 0& 0& 0 & 0& 0& 0 & 0 & 0& 0& 0\\
    -z^{3} & 0 & zx & -zx^{2} & xz^{2} &  0  &  -x^{2} &  x^{3} & 0& 0& 0 & 0 & 0& 0& 0&0&0& 0& 0& 0 & 0 & 0& y& 0& 0 & 0& 0& 0 & 0 & 0& 0& 0\\
     0 & -z^{3} & z^{2} & zy& 0 &  xz^{2}&  -xz & -xy & 0& 0& 0 & 0 & 0& 0& 0&0&0& 0& 0& 0 & 0 & 0& 0& y& 0 & 0& 0& 0 & 0 & 0& 0& 0\\
      0& 0& 0 & 0 & 0& 0& 0&0 & -xy & -x^{3} & -xy & 0 & -zy & -zx^{2} & -zy & 0&0& 0& 0& 0 & 0 & 0& 0& 0& y & 0& 0& 0 & 0 & 0& 0& 0\\
      0& 0& 0 & 0 & 0& 0& 0&0 & xz &   -x^{2}    &  0 & -xy & z^{2} & -zx & 0 & -zy &0& 0& 0& 0 & 0 & 0& 0& 0& 0 & y& 0& 0 & 0 & 0& 0& 0\\
      0& 0& 0 & 0 & 0& 0& 0&0 & xz^{2} &  0    &  -x^{2} &  x^{3} & z^{3} & 0 & -zx & zx^{2}&0& 0& 0& 0 & 0 & 0& 0& 0& 0 & 0& y& 0 & 0 & 0& 0& 0\\
       0& 0& 0 & 0 & 0& 0& 0&0 &0 &  xz^{2}&  -xz & -xy & 0 & z^{3} & -z^{2} & -zy&0& 0& 0& 0 & 0 & 0& 0& 0& 0 & 0& 0& y & 0 & 0& 0& 0\\
      0& 0& 0 & 0 & 0& 0& 0&0 & zy & zx^{2} & zy & 0 & -xy & -x^{3} & -xy & 0&0& 0& 0& 0 & 0 & 0& 0& 0& 0 & 0& 0& 0 & y & 0& 0& 0\\
     0& 0& 0 & 0 & 0& 0& 0&0 & -z^{2} & zx & 0 & zy & xz &   -x^{2}    &  0 & -xy&0& 0& 0& 0 & 0 & 0& 0& 0& 0 & 0& 0& 0 & 0 & y& 0& 0\\
      0& 0& 0 & 0 & 0& 0& 0&0 & -z^{3} & 0 & zx & -zx^{2} & xz^{2} &  0  &  -x^{2} & x^{3}&0& 0& 0& 0 & 0 & 0& 0& 0& 0 & 0& 0& 0 & 0 & 0& y& 0\\
      0& 0& 0 & 0 & 0& 0& 0&0 & 0 & -z^{3} & z^{2} & zy& 0 &  xz^{2}&  -xz & -xy &0& 0& 0& 0 & 0 & 0& 0& 0& 0 & 0& 0& 0 & 0 & 0& 0& y
 \end{pmatrix}
  $}
\end{gather*}

\begin{gather*}
  \setlength{\arraycolsep}{.60\arraycolsep}
  \renewcommand{\arraystretch}{1.5}
  \text{\footnotesize$\displaystyle
    \begin{pmatrix}
        y& 0& 0& 0 & 0 & 0& 0& 0& 0 & 0& 0& 0 & 0 & 0& 0& 0 & -x^{2} & x^{3} & xy & 0 & zx & -zx^{2} & -zy & 0 & 0& 0& 0 & 0 & 0& 0& 0&0 \\
    0& y& 0& 0 & 0 & 0& 0& 0& 0 & 0& 0& 0 & 0 & 0& 0& 0 &-xz &   -xy    &  0 & xy & z^{2} & zy & 0 & -zy & 0& 0& 0 & 0 & 0& 0& 0&0 \\
    0& 0& y& 0 & 0 & 0& 0& 0& 0 & 0& 0& 0 & 0 & 0& 0& 0 &-xz^{2} &  0 &  -xy &  -x^{3} & z^{3} & 0 & zy & zx^{2} & 0& 0& 0 & 0 & 0& 0& 0&0 \\
     0& 0& 0& y & 0 & 0& 0& 0& 0 & 0& 0& 0 & 0 & 0& 0& 0 &0 &  -xz^{2}&  xz & -x^{2} & 0 & z^{3} & -z^{2} & zx & 0& 0& 0 & 0 & 0& 0& 0&0 \\
     0& 0& 0& 0 & y & 0& 0& 0& 0 & 0& 0& 0 & 0 & 0& 0& 0 &-zx & zx^{2} & zy & 0 & -x^{2} & x^{3} & xy & 0 & 0& 0& 0 & 0 & 0& 0& 0&0 \\
     0& 0& 0& 0 & 0 & y& 0& 0& 0 & 0& 0& 0 & 0 & 0& 0& 0 &-z^{2} & -zy & 0 & zy & -xz & -xy & 0 & xy & 0& 0& 0 & 0 & 0& 0& 0&0 \\
    0& 0& 0& 0 & 0 & 0& y& 0& 0 & 0& 0& 0 & 0 & 0& 0& 0 & -z^{3} & 0 & -zy & -zx^{2} & -xz^{2} & 0 & -xy & -x^{3} & 0& 0& 0 & 0 & 0& 0& 0&0 \\
      0& 0& 0& 0 & 0 & 0& 0& y& 0 & 0& 0& 0 & 0 & 0& 0& 0 &0 & -z^{3} & z^{2} & -zx &  0 &  -xz^{2}&  xz & -x^{2} & 0& 0& 0 & 0 & 0& 0& 0&0 \\
      0& 0& 0& 0 & 0 & 0& 0& 0& y & 0& 0& 0 & 0 & 0& 0& 0 &0& 0& 0 & 0 & 0& 0& 0&0 & -x^{2} & x^{3} & xy & 0 & zx & -zx^{2} & -zy & 0\\
     0& 0& 0& 0 & 0 & 0& 0& 0& 0 & y& 0& 0 & 0 & 0& 0& 0 &0& 0& 0 & 0 & 0& 0& 0&0 & -xz &   -xy    &  0 & xy & z^{2} & zy & 0 & -zy \\
     0& 0& 0& 0 & 0 & 0& 0& 0& 0 & 0& y& 0 & 0 & 0& 0& 0 &0& 0& 0 & 0 & 0& 0& 0&0 & -xz^{2} &  0 &  -xy &  -x^{3} & z^{3} & 0 & zy & zx^{2}\\
    0& 0& 0& 0 & 0 & 0& 0& 0& 0 & 0& 0& y & 0 & 0& 0& 0 & 0& 0& 0 & 0 & 0& 0& 0&0 & 0 &  -xz^{2}&  xz & -x^{2} & 0 & z^{3} & -z^{2} & zx\\
     0& 0& 0& 0 & 0 & 0& 0& 0& 0 & 0& 0& 0 & y & 0& 0& 0 &0& 0& 0 & 0 & 0& 0& 0&0 & -zx & zx^{2} & zy & 0 & -x^{2} & x^{3} & xy & 0\\
     0& 0& 0& 0 & 0 & 0& 0& 0& 0 & 0& 0& 0 & 0 & y& 0& 0 &0& 0& 0 & 0 & 0& 0& 0&0 & -z^{2} & -zy & 0 & zy & -xz & -xy & 0 & xy \\
     0& 0& 0& 0 & 0 & 0& 0& 0& 0 & 0& 0& 0 & 0 & 0& y& 0 &0& 0& 0 & 0 & 0& 0& 0&0 & -z^{3} & 0 & -zy & -zx^{2} & -xz^{2} & 0 & -xy & -x^{3} \\
      0& 0& 0& 0 & 0 & 0& 0& 0& 0 & 0& 0& 0 & 0 & 0& 0& y &0& 0& 0 & 0 & 0& 0& 0&0 &0 & -z^{3} & z^{2} & -zx &  0 &  -xz^{2}&  xz & -x^{2}\\
       xy & x^{3} & xy & 0 & zy & zx^{2} & zy & 0 & 0& 0& 0 & 0 & 0& 0& 0&0 &x & 0& 0& 0 & 0 & 0& 0& 0& 0 & 0& 0& 0 & 0 & 0& 0& 0\\
    -xz &   x^{2}    &  0 & xy & -z^{2} & zx & 0 & zy & 0& 0& 0 & 0 & 0& 0& 0&0&0 & x& 0& 0 & 0 & 0& 0& 0& 0 & 0& 0& 0 & 0 & 0& 0& 0\\
   -xz^{2} &  0    &  x^{2} &  -x^{3} & -z^{3} & 0 & zx & -zx^{2} & 0& 0& 0 & 0 & 0& 0& 0&0&0 & 0& x& 0 & 0 & 0& 0& 0& 0 & 0& 0& 0 & 0 & 0& 0& 0\\
     0 &  -xz^{2}&  xz & xy & 0 & -z^{3} & z^{2} & zy & 0& 0& 0 & 0 & 0& 0& 0&0&0 & 0& 0& x & 0 & 0& 0& 0& 0 & 0& 0& 0 & 0 & 0& 0& 0\\
    -zy & -zx^{2} & -zy & 0 & xy & x^{3} & xy & 0 & 0& 0& 0 & 0 & 0& 0& 0&0&0 & 0& 0& 0 & x & 0& 0& 0& 0 & 0& 0& 0 & 0 & 0& 0& 0\\
    z^{2} & -zx & 0 & -zy & -xz &   x^{2}    &  0 & xy & 0& 0& 0 & 0 & 0& 0& 0&0&0 & 0& 0& 0 & 0 & x& 0& 0& 0 & 0& 0& 0 & 0 & 0& 0& 0\\
    z^{3} & 0 & -zx & zx^{2} & -xz^{2} &  0  &  x^{2} &  -x^{3} & 0& 0& 0 & 0 & 0& 0& 0&0&0 & 0& 0& 0 & 0 & 0& x& 0& 0 & 0& 0& 0 & 0 & 0& 0& 0\\
     0 & z^{3} & -z^{2} & -zy& 0 &  -xz^{2}&  xz & xy & 0& 0& 0 & 0 & 0& 0& 0&0&0 & 0& 0& 0 & 0 & 0& 0& x& 0 & 0& 0& 0 & 0 & 0& 0& 0\\
      0& 0& 0 & 0 & 0& 0& 0&0 & xy & x^{3} & xy & 0 & zy & zx^{2} & zy & 0&0 & 0& 0& 0 & 0 & 0& 0& 0& x & 0& 0& 0 & 0 & 0& 0& 0\\
      0& 0& 0 & 0 & 0& 0& 0&0 & -xz &   x^{2}    &  0 & xy & -z^{2} & zx & 0 & zy &0 & 0& 0& 0 & 0 & 0& 0& 0& 0 & x& 0& 0 & 0 & 0& 0& 0\\
      0& 0& 0 & 0 & 0& 0& 0&0 &-xz^{2} &  0    &  x^{2} &  -x^{3} & -z^{3} & 0 & zx & -zx^{2}&0 & 0& 0& 0 & 0 & 0& 0& 0& 0 & 0& x& 0 & 0 & 0& 0& 0\\
       0& 0& 0 & 0 & 0& 0& 0&0 &0 &  -xz^{2}&  xz & xy & 0 & -z^{3} & z^{2} & zy&0 & 0& 0& 0 & 0 & 0& 0& 0& 0 & 0& 0& x & 0 & 0& 0& 0\\
      0& 0& 0 & 0 & 0& 0& 0&0 &-zy & -zx^{2} & -zy & 0 & xy & x^{3} & xy & 0&0 & 0& 0& 0 & 0 & 0& 0& 0& 0 & 0& 0& 0 & x & 0& 0& 0\\
     0& 0& 0 & 0 & 0& 0& 0&0 & z^{2} & -zx & 0 & -zy & -xz &   x^{2}    &  0 & xy&0 & 0& 0& 0 & 0 & 0& 0& 0& 0 & 0& 0& 0 & 0 & x& 0& 0\\
      0& 0& 0 & 0 & 0& 0& 0&0 &z^{3} & 0 & -zx & zx^{2} & -xz^{2} &  0  &  x^{2} &  -x^{3}&0 & 0& 0& 0 & 0 & 0& 0& 0& 0 & 0& 0& 0 & 0 & 0& x& 0\\
      0& 0& 0 & 0 & 0& 0& 0&0 & 0 & z^{3} & -z^{2} & -zy& 0 &  -xz^{2}&  xz & xy&0 & 0& 0& 0 & 0 & 0& 0& 0& 0 & 0& 0& 0 & 0 & 0& 0& x
 \end{pmatrix}
  $}
\end{gather*}

\end{enumerate}     
    
\end{example}

\begin{corollary} \label{Cor. improved algo on t+gh}
  Let $f$ be a simple summand-reducible polynomial with $n+1$ monomials in its expanded form. If $f$ can be written as $t + gh$ where $t$ is a monomial and $gh$ is a product of sums of monomials, where $p$ is the number of monomials in $g$ and $q$ is the number of monomials in $h$, then the improved algorithm of theorem \ref{thm improved standard method} produces a matrix factorization of $f$ which is $2^{n-(p+q)}=2^{pq -(p+q)}$ times smaller in size than what one would obtain using the standard method.
\end{corollary}

\begin{proof}
 Let $f$ be a simple summand-reducible polynomial with $n+1$ monomials in its expanded form and such that $f=t+gh$ as in the assumption of the theorem. Let $p$ be the number of monomials in $g$ and $q$ the number of monomials in $h$, then $n=pq$ and the standard method would produce a factorization of $f$ of size $2^{(n+1)-1}=2^{(pq+1) - 1}=2^{pq}$. \\Let us now use the improved algorithm to find the size of matrix factors of $f$. Since $f=t+gh$, we use the standard method to find matrix factorizations of $g$ and $h$ and find respectively  factorizations of size $2^{p - 1}$  and $2^{q - 1}$. Next, we use the multiplicative tensor product of matrix factorizations $(\widetilde{\otimes})$ to find a factorization of size $2(2^{p - 1})(2^{q - 1})=2^{p + q - 1}$ for the product $gh$.
 Thereafter, use the tensor product of matrix factorizations $(\widehat{\otimes})$ to find a matrix factorization of $f$ of size $2(1)(2^{p + q - 1})=2^{p + q}$. \\
 Now we compare the sizes obtained by the two methods. WLOG, $p\geq 2$ and $q\geq 3$, so $pq > p + q$ and $2^{pq} > 2^{p + q }$.\\
 Finally, $2^{pq} \div 2^{p + q}= 2^{(pq)- (p + q)}$ as desired.
\end{proof}
This theorem actually shows that the higher the product $pq$, the smaller the size of the factorization we obtain as compared to what one would obtain using the standard method.\\
We now generalize definition \ref{defn simple summand-reducible polynomial} and define the class of \textit{summand-reducible} polynomials which is made up of polynomials in which some monomials can be factorized in a nice way, hence allowing the polynomial to be written with less summands.

\begin{definition}\label{defn summand reducible polynomials}
  A polynomial $f$ is said to be \textbf{summand-reducible} if it can be written in the form $$f=t_{1}+\cdots + t_{s}+ g_{11}\cdots g_{1m_{1}} + \cdots + g_{l1}\cdots g_{lm_{l}},$$ where:
  \begin{enumerate}
    \item $\bullet$ If $s=0$, then there exist at least two products $g_{11}\cdots g_{1m_{1}}$ and $g_{21}\cdots g_{2m_{2}}$ in $f$.\\
        $\bullet$ If $s\neq 0$, then there exists at least one product $g_{11}\cdots g_{1m_{1}}$ in $f$.
    \item For $i=1,\cdots,s$; each $t_{i}$ is a monomial and so $t_{i}=h_{i1}h_{i2}$, where $h_{i1}$ and $h_{i2}$ are products of variables possibly raised to some power.
    \item For $j=1,\cdots, l$; each $g_{j1}\cdots g_{jm_{j}}$ is a product of sums of monomials, such that if it is expanded, $g_{j1}\cdots g_{jm_{j}}$ would have more monomials than the number that appears in the factor form $g_{j1}\cdots g_{jm_{j}}$.
    \item For $1\leq j \leq l$, at least one of the products $g_{j1}\cdots g_{jm_{j}}$ has at least two factors.
  \end{enumerate}
\end{definition}
\begin{definition}
  A polynomial $f$ is said to be \textbf{summand-reduced} if it is in the form $f=t_{1}+\cdots + t_{s}+ g_{11}\cdots g_{1m_{1}} + \cdots + g_{l1}\cdots g_{lm_{l}}$ described in definition \ref{defn summand reducible polynomials}.
\end{definition}

\begin{example}
   $f= xy +(xy+x^{2}z+yz^{2})(x^{2}+y^{2}) +(yz + xy^{2} +x^{2})(x^{3}z^{2}+ yx + y^{2})$ is a summand-reduced polynomial.
  \end{example}
\begin{example}
  Every simple summand-reducible polynomial is a summand-reducible polynomial. So example \ref{exples of simple summand-red polyn} also gives examples of summand-reducible polynomials.\\
\end{example}
We now generalize theorem \ref{thm improved standard method}.
\begin{theorem} \label{thm improved algo for summand-red polyn}
Let $f=t_{1}+\cdots + t_{s}+ g_{11}\cdots g_{1m_{1}} + \cdots + g_{l1}\cdots g_{lm_{l}}$ be a summand-reducible polynomial. Let $p_{ji}$ be the number of monomials in $g_{ji}$. Then
there is an improved version of the standard method for factoring $f$ which produces factorizations of size $$2^{\prod_{i=1}^{m_{1}}p_{1i} + \cdots + \prod_{i=1}^{m_{l}}p_{li} - (\sum_{i=1}^{m_{1}}p_{1i} + \cdots + \sum_{i=1}^{m_{l}}p_{li})}$$ times smaller than the size one would normally obtain with the standard method.
\end{theorem}

\begin{proof}
  First, we construct the algorithm, then we prove that the resulting matrix factorizations (for summand-reducible polynomials) are $$2^{\prod_{i=1}^{m_{1}}p_{1i} + \cdots + \prod_{i=1}^{m_{l}}p_{li} - (\sum_{i=1}^{m_{1}}p_{1i} + \cdots + \sum_{i=1}^{m_{l}}p_{li})}$$ times smaller in size than what one would obtain with the standard method.\\
  We inductively construct the matrix factorizations of summand-reduced polynomials using the tensor products $\widehat{\otimes}$ and $\widetilde{\otimes}$  that were not existing in the 1980s when the standard method was constructed.\\
  The algorithm we propose here is just an improvement of the one we gave for simple summand-reducible polynomials (cf. proof of theorem \ref{thm improved standard method}).

Let $f=t_{1}+\cdots + t_{s}+ g_{11}\cdots g_{1m_{1}} + \cdots + g_{l1}\cdots g_{lm_{l}}$ be a summand-reducible polynomial. Let $p_{ji}$ be the number of monomials in $g_{ji}$. \\
If $\forall k\in \{1,\cdots,s\}$, $t_{k}= 0$, then do:
\begin{enumerate}
  \item For each $j\in \{1,\cdots, l\}$ and $i\in \{1,\cdots, m_{j}\}$, use the standard method to find a matrix factorization of $g_{ji}$ of size $2^{p_{ji}-1}$.
  \item Next, for each $j\in \{1,\cdots, l\}$; use the multiplicative tensor product of matrix factorizations $\widetilde{\otimes}$ (or its variant) to find a matrix factorization of $g_{j1}\cdots g_{jm_{j}}$ of size $$(2^{m_{j}-1})(2^{\sum_{i=1}^{m_{j}}p_{ji}-m_{j}})= 2^{\sum_{i=1}^{m_{j}}p_{ji}-1}$$
  \item Now use the tensor product of matrix factorizations $\widehat{\otimes}$ to find a matrix factorization of $g_{11}\cdots g_{1m_{1}} + \cdots + g_{l1}\cdots g_{lm_{l}}$ of size
      $$(2^{l-1})(\prod_{j=1}^{l}2^{\sum_{i=1}^{m_{j}}p_{ji}-1}) = 2^{l-1+ \sum_{i=1}^{m_{1}}p_{1i} + \cdots + \sum_{i=1}^{m_{l}}p_{li} - l}= 2^{ \sum_{i=1}^{m_{1}}p_{1i} + \cdots + \sum_{i=1}^{m_{l}}p_{li} - 1}.$$
      Let us find the size of matrices the standard method would produce for $$g_{11}\cdots g_{1m_{1}} + \cdots + g_{l1}\cdots g_{lm_{l}}.$$
      Let $n_{j}=$number of monomials in the expanded form of the $j^{th}$ product $g_{j1}\cdots g_{jm_{j}}$. Then $n_{j}=\prod_{i=1}^{m_{j}} p_{ji}$. Hence, the number of monomials in the expanded form of $g_{11}\cdots g_{1m_{1}} + \cdots + g_{l1}\cdots g_{lm_{l}}$ would be $\sum_{j=1}^{l} n_{j}= \sum_{j=1}^{l} \prod_{i=1}^{m_{j}} p_{ji}$. \\So the size of factorizations produced by the standard method would be $2^{(\sum_{j=1}^{l} \prod_{i=1}^{m_{j}}p_{ji})-1}.$ \\
      Thus, the factorizations produced by our improved algorithm would be
       $$2^{(\sum_{j=1}^{l} \prod_{i=1}^{m_{j}}p_{ji})-1}\div 2^{ \sum_{i=1}^{m_{1}}p_{1i} + \cdots + \sum_{i=1}^{m_{l}}p_{li} - 1}= 2^{(\sum_{j=1}^{l} \prod_{i=1}^{m_{j}}p_{ji})- (\sum_{i=1}^{m_{1}}p_{1i} + \cdots + \sum_{i=1}^{m_{l}}p_{li})}$$
      times smaller in size than the factorizations produced by the standard method.
  \item If there exists $k\in \{1,\cdots,s\}$ such that $t_{k}\neq 0$, then use the standard method to inductively find a matrix factorization $(A,B)$ of $t_{1}+\cdots + t_{s}$ of size $2^{s-1}$.
  \item Then do steps 1), 2) and 3) above to find a matrix factorization $(C,D)$ of $g_{11}\cdots g_{1m_{1}} + \cdots + g_{l1}\cdots g_{lm_{l}}$ of size $2^{ \sum_{i=1}^{m_{1}}p_{1i} + \cdots + \sum_{i=1}^{m_{l}}p_{li} - 1}$.
  \item Now, use $\widehat{\otimes}$ to find a matrix factorization $(A,B) \widehat{\otimes} (C,D)$ of $f=t_{1}+\cdots + t_{s}+ g_{11}\cdots g_{1m_{1}} + \cdots + g_{l1}\cdots g_{lm_{l}}$ of size $$2(2^{s-1})(2^{\sum_{i=1}^{m_{1}}p_{1i} + \cdots + \sum_{i=1}^{m_{l}}p_{li} - 1})= 2^{\sum_{i=1}^{m_{1}}p_{1i} + \cdots + \sum_{i=1}^{m_{l}}p_{li} + s - 1}.$$
      Note that $f$ in expanded form has $$\sum_{j=1}^{l} n_{j}+ s= (\sum_{j=1}^{l} \prod_{i=1}^{m_{j}}p_{ji}) + s$$ monomials and so the standard method would produce factorizations of size $2^{(\sum_{j=1}^{l} \prod_{i=1}^{m_{j}}p_{ji}) + s - 1}$.\\ Hence the factorizations our improved algorithm produce are $$2^{(\sum_{j=1}^{l} \prod_{i=1}^{m_{j}}p_{ji}) + s - 1}\div 2^{ \sum_{i=1}^{m_{1}}p_{1i} + \cdots + \sum_{i=1}^{m_{l}}p_{li} + s - 1}= 2^{(\sum_{j=1}^{l} \prod_{i=1}^{m_{j}}p_{ji})- (\sum_{i=1}^{m_{1}}p_{1i} + \cdots + \sum_{i=1}^{m_{l}}p_{li})}$$
      times smaller in size than the factorizations produced by the standard method. QED.

\end{enumerate}
\end{proof}
\begin{example} \label{exple: first one after the main thm}
  Let $f= xy +(xy+x^{2}z+yz^{2})(x^{2}+z^{2}) +(yz + xy^{2} +x^{2})(x^{3}z^{2}+ yx + y^{2})$.\\
  $f$ in expanded form has $1 + 3\times 2 + 3\times 3=16$ monomials and so the standard method will produce factorizations of size $2^{16-1}=2^{15}$.\\
  The foregoing theorem gives us a way to find the size of the factorizations we would obtain using the improved algorithm.\\
  In fact, for this example: $s=1$, $l=2$, $m_{1}=2$, $m_{2}=2$ $p_{11}=p_{21}=p_{22}=3$ and $p_{12}=2$. So, our algorithm would produce factorizations of size $2^{p_{11}+p_{12}+p_{21}+p_{22}+s-1}=2^{3+2+3+3+1-1}=2^{11}.$\\
  Hence, from the theorem we deduce that the improved algorithm produces factorizations of size $2^{15} \div 2^{11}=2^{4}=16$ times smaller than what the standard method produces!
\end{example}

\begin{example} \label{exple: matrix factors presented, part II}
Use the improved algorithm to factorize the polynomial $f(x)= x^{3}y^{2} + (xy+x^{2}z+yz^{2})(xz+y^{2}+y^{2}z)$ and compare the size of the matrix factors with the one obtained using the standard method.\\
Since $f$ in its expanded form has $1+3\times 3=1+9=10$ monomials, the size of matrix factors obtained using the standard method would be $2^{10-1}=2^{9}=512$. We can use Theorem \ref{thm improved algo for summand-red polyn} as we did in Example \ref{exple: first one after the main thm} to find that the size of matrix factors of $f$ using the improved algorithm is $2^{3+3+1-1}=2^{6}=64$, that is $\frac{512}{64}=8$ times smaller than the size obtained using the standard method.\\
In the sequel, we are going to use this algorithm to find matrix factors of $f$ and we will see that they are of size $64$.\\
 Let $g=xy+x^{2}z+yz^{2}$ and $s=xz+y^{2}+y^{2}z$, so that $f=x^{3}y^{2} + gs$.\\

 In Example \ref{exple: good matrix facto of g}, we used the standard method to find a matrix factorization of the polynomial $g$:

\[
  P=(\phi_{g},\psi_{g})=(\begin{bmatrix}
    x & -x^{2} & -y & 0\\
    z &   y    &  0 & -y\\
 z^{2} &  0    &  y &  x^{2}\\
     0 &  z^{2}&  -z & x
  \end{bmatrix},\begin{bmatrix}
     y & x^{2} &  y & 0\\
    -z &   x   &  0 & y\\
 -z^{2} &  0    &  x &  -x^{2}\\
     0 &  -z^{2}&  z & y
  \end{bmatrix})
\]

Let us find a matrix factorization of $s=xz+y^{2}+y^{2}z=d+y^{2}z$, where $d=xz+y^{2}$. Using the standard method, we find that
\[
  (\begin{bmatrix}
    x & -y \\
    y & z
  \end{bmatrix},\begin{bmatrix}
    z & y \\
    -y & x
  \end{bmatrix})
\]
  is a matrix factorization of $d=xz+y^{2}$. Thus, using the standard method, a matrix factorization of the polynomial $s$ is

\[
  N=(\phi_{s},\psi_{s})=(\begin{bmatrix}
    x & -y & -y^{2} & 0\\
    y &  z  &  0 & -y^{2}\\
   z &  0  &  z &  y\\
     0 &  z&  -y & x
  \end{bmatrix},\begin{bmatrix}
     z & y &  y^{2} & 0\\
    -y &   x   &  0 & y^{2}\\
 -z &  0    &  x &  -y\\
     0 &  -z &  y & z
  \end{bmatrix})
\] 
  
According to the proof of Theorem \ref{thm improved algo for summand-red polyn} to find a matrix factorization for $f$, we need to:
\begin{enumerate}
  \item First of all find a matrix factorization of the product $gs$ using the multiplicative  tensor product $\widetilde{\otimes}$. By Lemma \ref{lemma size of X tensor Y}, the matrix factors of the product $gs$ will be of size $2(4)(4)=32$ since $P$ and $N$ (which are respectively matrix factorizations of $g$ and $s$) are each of size $4$. \\
       We proceed as follows:\\
$P \widetilde{\otimes} N = (\phi_{g},\psi_{g})\widetilde{\otimes}(\phi_{s},\psi_{s})=(\phi_{gs},\psi_{gs})$ where $$(\phi_{gs},\psi_{gs})= (\begin{bmatrix}
    \phi_{g}\otimes \phi_{s} & 0 \\
    0             &   \phi_{g}\otimes \phi_{s}
  \end{bmatrix},\begin{bmatrix}
   \psi_{g}\otimes \psi_{s} & 0 \\
    0             &   \psi_{g}\otimes \psi_{s}
  \end{bmatrix})$$
with \\
$\,\,\,\,\,\,\,\,\,\,\phi_{g}\otimes \phi_{s}\,\,\,\,\,\,\,\,\,\,\,\,\,\,\,\,\,\,=
\begin{bmatrix}
  x & -x^{2} & -y & 0\\
  z  &  y  &  0 & -y\\
 z^{2} &  0  &  y &  x^{2}\\
     0 &  z^{2}&  -z & x \end{bmatrix} \otimes \begin{bmatrix}
      x & -y & -y^{2} & 0\\
    y &  z  &  0 & -y^{2}\\
   z &  0  &  z &  y\\
     0 &  z&  -y & x
  \end{bmatrix}$

  \begin{gather*}
  \setlength{\arraycolsep}{1.0\arraycolsep}
  \renewcommand{\arraystretch}{1.5}
  \text{\footnotesize$\displaystyle
    i.e., \,\phi_{g}\otimes \phi_{s}=\begin{pmatrix}
      x^{2} & -xy & -xy^{2} & 0 & -x^{3} & x^{2}y & x^{2}y^{2} & 0 & -yx & y^{2} & y^{3} & 0 & 0& 0& 0&0 \\
    xy &   xz    &  0 & -xy^{2} & -x^{2}y & -x^{2}z & 0 & x^{2}y^{2} & -y^{2}& -yz & 0 & y^{3} & 0& 0& 0&0 \\
    xz &  0 &  xz &  xy & -x^{2}z & 0 & -x^{2}z & -x^{2}y & -yz& 0& -yz & -y^{2} & 0& 0& 0&0 \\
     0 &  xz &  -xy & x^{2} & 0 & -x^{2}z & x^{2}y & -x^{3} & 0& -yz& y^{2} & -yx & 0& 0& 0&0 \\
     zx & -zy & -zy^{2} & 0 & yx & -y^{2} & -y^{3} & 0 & 0& 0& 0 & 0 & -yx& y^{2}& y^{3}&0 \\
     zy & z^{2} & 0 & -zy^{2} & y^{2} & yz & 0 & -y^{3} & 0& 0& 0 & 0 & -y^{2}& -yz& 0&y^{3} \\
     z^{2} & 0 & z^{2} & zy & yz & 0 & yz & y^{2} & 0& 0& 0 & 0 & -yz& 0& -yz&-y^{2} \\
      0 & z^{2} & -zy & zx &  0 &  yz &  -y^{2} & yx & 0 & 0& 0& 0&0 & -yz& y^{2} & -yx\\
      z^{2}x & -z^{2}y& -z^{2}y^{2} & 0 & 0& 0& 0&0 & yx & -y^{2} & -y^{3} & 0 & x^{3} & -x^{2}y & -x^{2}y^{2} & 0\\
     z^{2}y& z^{3}& 0 & -z^{2}y^{2} & 0& 0& 0&0 & y^{2} &   yz    &  0 & -y^{3} & x^{2}y & x^{2}z & 0 & -x^{2}y^{2} \\
     z^{3}& 0& z^{3} & z^{2}y & 0& 0& 0&0 & yz &  0 &  yz &  y^{2} & x^{2}z & 0 & x^{2}z & x^{2}y\\
     0& z^{3}& -z^{2}y & z^{2}x & 0& 0& 0&0 & 0 &  yz &  -y^{2} & yx & 0 & x^{2}z & -x^{2}y & x^{3}\\
     0& 0& 0 & 0 & z^{2}x & -z^{2}y& -z^{2}y^{2}&0 & -zx & zy & zy^{2} & 0 & x^{2} & -xy & -xy^{2} & 0\\
     0& 0& 0 & 0 & z^{2}y& z^{3}& 0&-z^{2}y^{2} &-zy  & -z^{2} & 0 & zy^{2} & xy & xz & 0 & -xy^{2} \\
     0& 0& 0 & 0 & z^{3}& 0& z^{3}&z^{2}y & -z^{2} & 0 & -z^{2} & -zy & xz & 0 & xz & xy \\
      0& 0& 0 & 0 & 0& z^{3}& -z^{2}y& z^{2}x & 0& -z^{2} & zy & -zx &  0 &  xz&  -xy & x^{2}
 \end{pmatrix}
  $}
\end{gather*}
  
  And 
  $\,\,\,\,\,\,\,\,\,\,\psi_{g}\otimes \psi_{s}\,\,\,\,\,\,\,\,\,\,\,\,\,\,\,\,\,\,=
  (\begin{bmatrix}
     y & x^{2} &  y & 0\\
    -z &   x   &  0 & y\\
 -z^{2} &  0    &  x &  -x^{2}\\
     0 &  -z^{2}&  z & y
  \end{bmatrix}\otimes \begin{bmatrix}
     z & y &  y^{2} & 0\\
    -y &   x   &  0 & y^{2}\\
 -z &  0    &  x &  -y\\
     0 &  -z &  y & z
  \end{bmatrix})$
  
   \begin{gather*}
  \setlength{\arraycolsep}{1.0\arraycolsep}
  \renewcommand{\arraystretch}{1.5}
  \text{\footnotesize$\displaystyle
    i.e., \,\psi_{g}\otimes \psi_{s}=\begin{pmatrix}
       yz& y^{2} & y^{3} & 0 & x^{2}z & x^{2}y & x^{2}y^{2} & 0 & yz & y^{2} & y^{3} & 0 & 0& 0& 0&0 \\
    -y^{2} &   yx   &  0 & y^{3} & -x^{2}y & x^{3} & 0 & x^{2}y^{2} & -y^{2}& yx & 0 & y^{3} & 0& 0& 0&0 \\
    -yz &  0 &  yx &  -y^{2} & -x^{2}z & 0 & x^{3} & -x^{2}y & -yz& 0& yx & -y^{2} & 0& 0& 0&0 \\
     0 &  -yz &  y^{2} & yz  & 0 & -x^{2}z & x^{2}y & x^{2}z & 0& -yz& y^{2} & yz & 0& 0& 0&0 \\
     -z^{2} & -zy & -zy^{2} & 0 & xz & xy & xy^{2} & 0 & 0& 0& 0 & 0 & yz& y^{2}& y^{3}&0 \\
     zy & -zx & 0 & -zy^{2} & -xy & x^{2} & 0 & xy^{2} & 0& 0& 0 & 0 & -y^{2}& yx& 0&y^{3} \\
     z^{2} & 0 & -zx & zy & -xz & 0 & x^{2} & -xy & 0& 0& 0 & 0 & -yz& 0& yx &-y^{2} \\
      0 & z^{2} & -zy & -z^{2} &  0 &  -xz &  xy & xz & 0 & 0& 0& 0&0 & -yz& y^{2} & yz\\
      -z^{3} & -z^{2}y& -z^{2}y^{2} & 0 & 0& 0& 0&0 & xz & xy & xy^{2} & 0 & -x^{2}z & -x^{2}y & -x^{2}y^{2} & 0\\
     z^{2}y& -z^{2}x& 0 & -z^{2}y^{2} & 0& 0& 0&0 & -xy &   x^{2}    &  0 & xy^{2} & x^{2}y & -x^{3} & 0 & -x^{2}y^{2} \\
     z^{3}& 0& -z^{2}x & z^{2}y & 0& 0& 0&0 & -xz &  0 &  x^{2} &  -xy & x^{2}z & 0 & -x^{3} & x^{2}y\\
     0& z^{3}& -z^{2}y & -z^{3} & 0& 0& 0&0 & 0 &  -xz &  xy & xz & 0 & x^{2}z & -x^{2}y & -x^{2}z\\
     0& 0& 0 & 0 & -z^{3} & -z^{2}y& -z^{2}y^{2}&0 & z^{2} & zy & zy^{2} & 0 & yz & y^{2} & y^{3} & 0\\
     0& 0& 0 & 0 & z^{2}y& -z^{2}x& 0&-z^{2}y^{2} &-zy  & zx & 0 & zy^{2} & -y^{2} & yx & 0 & y^{3} \\
     0& 0& 0 & 0 & z^{3}& 0& -z^{2}x &z^{2}y & -z^{2} & 0 & zx & -zy & -yz & 0 & yx & -y^{2} \\
      0& 0& 0 & 0 & 0& z^{3}& -z^{2}y& -z^{3} & 0& -z^{2} & zy & z^{2} &  0 &  -yz &  y^{2} & yz
 \end{pmatrix}
  $}
\end{gather*}

 \item Next, from the algorithm given in Theorem \ref{thm improved algo for summand-red polyn}, we now need to find a matrix factorization of $r=x^{3}y^{2}$ (which is the first summand in $f$). Evidently, $L=(\phi_{r},\psi_{r})=([x^{3}],[y^{2}])$ is a $1\times 1$ matrix factorization of $x^{3}y^{2}$.
\item Finally, from our algorithm we find a matrix factorization of $f$ by computing $L\widehat{\otimes} (P\widetilde{\otimes} N)$ which will be of size $2(1)(32)=64$ by Lemma 2.1 of \cite{fomatati2019multiplicative} since $L$ is of size $1$ and $(P\widetilde{\otimes} N)$ is of size $32$.\\
    We have:
    \begin{align*}
    L\widehat{\otimes} (P\widetilde{\otimes} N) &= (\phi_{r},\psi_{r})\widehat{\otimes} (\phi_{gs}, \psi_{gs})\\
&=( \begin{bmatrix}
    \phi_{r}\otimes 1_{32}  &  1_{1}\otimes \phi_{gs}      \\
   -1_{1}\otimes \psi_{gs}  &  \psi_{r}\otimes 1_{32}
\end{bmatrix},
\begin{bmatrix}
    \psi_{r}\otimes 1_{32}  &  -1_{1}\otimes \phi_{gs}     \\
    1_{1}\otimes \psi_{gs}  &  \phi_{r}\otimes 1_{32}
\end{bmatrix}
)\\
&=( \begin{bmatrix}
    x^{3}\otimes 1_{32}  &  1\otimes \phi_{gs}      \\
   -1\otimes \psi_{gs}  &  y^{2}\otimes 1_{32}
\end{bmatrix},
\begin{bmatrix}
    y^{2}\otimes 1_{32}  &  -1\otimes \phi_{gs}     \\
    1\otimes \psi_{gs}  &  x^{3}\otimes 1_{32}
\end{bmatrix}
)\\
&=( \begin{bmatrix}
    x^{3}\otimes 1_{32}  &   \phi_{gs}      \\
   - \psi_{gs}  &  y^{2}\otimes 1_{32}
\end{bmatrix},
\begin{bmatrix}
    y^{2}\otimes 1_{32}  &  -\phi_{gs}     \\
    \psi_{gs}  &  x^{3}\otimes 1_{32}
\end{bmatrix}
)\\
&=(\phi_{rgs},\psi_{rgs})
\end{align*}

Where:\\
$\bullet$ $x^{3}\otimes 1_{32}$ (respectively $y^{2}\otimes 1_{32}$) is a $32\times 32$ diagonal matrix with $x^{3}$ (respectively $y^{2}$) on its entire diagonal. \\
$\bullet$ $\phi_{gs}$ and $\psi_{gs}$ were computed above. \\
So, $(\phi_{rgs},\psi_{rgs})$ have been completely specified.
Hence, we found a $64\times 64$ matrix factorization of $f$ viz. a matrix factorization of $f$ of size 64.

\end{enumerate}

\end{example}

\section{Further research directions}
This research generated a number of other questions that we think are interesting.\\
\begin{enumerate}
\item Yoshino \cite{yoshino1998tensor} defined the tensor product ($\widehat{\otimes}$) of matrix factorizations and proved that if $X$ and $Y$ are respectively the matrix factorizations of power series $f$ and $g$, such that $X$ and $Y$ are indecomposable, it might still happen that $X\widehat{\otimes}Y$ be decomposable. He proved theorems (section 3 of \cite{yoshino1998tensor}) which give the bounds for the number of indecomposable components in the direct decomposition of $X\widehat{\otimes}Y$. It would be interesting to know if similar results could be derived for the tensor products proposed in this paper.
\item In this paper, we learned about factorization of polynomials using two $n\times n$ matrices ($n\geq 2$), it could be interesting to study how polynomials could be factorized using $m$ matrices of size $n$ where ($n\geq 2$) and $m\geq 3$. And one could also be interested in finding conditions under which a given polynomial has at most $r$ matrix factors, where $r\geq 2$.
\item So far, we have tensor products $\widehat{\otimes}$ and $\widetilde{\otimes}$ which have the ability to produce matrix factorizations of the sum $f+g$ and the product $fg$ respectively; from the matrix factorizations of the polynomials $f$ and $g$. An interesting question would be to find a tensor product which can produce a matrix factorization of the composition $f\circ g$ whenever this composition makes sense.
\end{enumerate}
\begin{quote}
  \textbf{Acknowledgments}
\end{quote}
Part of this work was carried out during my Ph.D. studies in mathematics at the University of Ottawa in Canada.
I am grateful to Prof. Dr. Richard Blute who was my Ph.D. supervisor for all the fruitful interactions. It is also a pleasure to thank the anonymous referee. \\ Finally, I gratefully acknowledge the financial support of the Queen Elizabeth Diamond Jubilee scholarship during my Ph.D. studies.



\bibliography{fomatati_ref}
\addcontentsline{toc}{section}
{References}
\end{document}